\documentclass[final]{elsarticle}
\usepackage[english]{babel}
\usepackage[utf8]{inputenc}
\usepackage{bm}
\usepackage{amsmath,amsthm,amsfonts}
\usepackage{graphicx}
\usepackage{color}
\usepackage{multirow}
\usepackage{subcaption}
\usepackage{soul}
\usepackage{mathtools}

\usepackage{algorithm2e}
\SetKwComment{Comment}{/* }{ */}
\RestyleAlgo{ruled}

\usepackage{xcolor}

\graphicspath{ {./figs/} }

\usepackage[mathlines]{lineno}
\newcommand*\patchAmsMathEnvironmentForLineno[1]{%
  \expandafter\let\csname old#1\expandafter\endcsname\csname #1\endcsname
  \expandafter\let\csname oldend#1\expandafter\endcsname\csname end#1\endcsname
  \renewenvironment{#1}%
     {\linenomath\csname old#1\endcsname}%
     {\csname oldend#1\endcsname\endlinenomath}}%
\newcommand*\patchBothAmsMathEnvironmentsForLineno[1]{%
  \patchAmsMathEnvironmentForLineno{#1}%
  \patchAmsMathEnvironmentForLineno{#1*}}%
\AtBeginDocument{%
\patchBothAmsMathEnvironmentsForLineno{equation}%
\patchBothAmsMathEnvironmentsForLineno{align}%
\patchBothAmsMathEnvironmentsForLineno{flalign}%
\patchBothAmsMathEnvironmentsForLineno{alignat}%
\patchBothAmsMathEnvironmentsForLineno{gather}%
\patchBothAmsMathEnvironmentsForLineno{multline}%
}

\begin{document}

\begin{frontmatter}

\title{Multiscale approximation and two-grid preconditioner for extremely anisotropic heat flow}
\author[tamu]{Maria Vasilyeva }
\address[tamu]{Department of Mathematics and Statistics, Texas A\&M University  - Corpus Christi, Corpus Christi, Texas, USA.}
\cortext[mycorrespondingauthor]{Corresponding author}
\ead{maria.vasilyeva@tamucc.edu}
\author[LANL]{Golo Wimmer}
\ead{gwimmer@lanl.gov}
\author[LANL]{Ben S. Southworth}
\ead{southworth@lanl.gov}
\address[LANL]{Theoretical Division, Los Alamos National Laboratory, P.O. Box 1663, Los Alamos, NM 87545 USA.}



\begin{abstract}
We consider anisotropic heat flow with extreme anisotropy, as arises in magnetized plasmas for fusion applications. Such problems pose significant challenges in both obtaining an accurate approximation as well in the construction of an efficient solver. In both cases, the underlying difficulty is 
in forming an accurate approximation of temperature fields that follow the direction of complex, non-grid-aligned magnetic fields. In this work, we construct a highly accurate coarse grid approximation using spectral multiscale basis functions based on local anisotropic normalized Laplacians. We show that the local generalized spectral problems yield local modes that align with magnetic fields, and provide an excellent coarse-grid approximation of the problem. We then utilize this spectral coarse space as an approximation in itself, and as the coarse-grid in a two-level spectral preconditioner. Numerical results are presented for several magnetic field distributions and anisotropy ratios up to $10^{12}$, showing highly accurate results with a large system size reduction, and two-grid preconditioning that converges in $\mathcal{O}(1)$ iterations, independent of anisotropy. 
\end{abstract}



\end{frontmatter}

\section{Introduction}
 
In magnetic confinement fusion, heat is transported along magnetic field lines up to 10 orders of magnitude faster than across \cite{gunter2007finite}. Together with typical time step sizes considered for magnetohydrodynamic (MHD) simulations, the corresponding extremely anisotropic diffusion equation poses numerical challenges both with respect to accuracy and solver efficiency. This holds true especially for scenarios in which the magnetic field topology is not aligned with the numerical mesh, which will generally be the case for more complex MHD simulations including magnetic islands and stochastic field configurations. In particular, low accuracy discretizations are typically characterized by numerical cross-diffusion of parallel heat flux in the perpendicular direction, leading to spuriously low confinement times.


To avoid this, higher-order approximations are often considered, e.g. \cite{sovinec2004nonlinear}. A high-order interior penalty discontinuous Galerkin finite element scheme is presented in \cite{green2022efficient}, and a hybridized discontinuous Galerkin scheme is considered in \cite{giorgiani2020high}. In \cite{gunter2007finite}, a finite element analog to the second-order, finite difference scheme for the solution of the anisotropic diffusion equation in strongly magnetized plasmas is presented with a fourth-order extension. Other approaches include mixed Virtual Element Methods as considered in \cite{berrone2023mixed}, and nonlinear constrained finite element approximations in \cite{kuzmin2009constrained}. Due to the highly directional nature of problems with extreme anisotropy ratios, as well as the effective decoupling of solutions on different field lines, techniques from discretizing hyperbolic equations have also proven useful. In \cite{sharma2007preserving}, the authors show that standard algorithms for anisotropic diffusion based on centered differencing do not preserve monotonicity and can lead to the violation of the entropy constraints of the second law of thermodynamics. To solve these problems, algorithms based on slope limiters, analogous to those used in second-order schemes for hyperbolic equations, are proposed. A high-order finite difference solver for anisotropic diffusion problems based on the first-order hyperbolic system method is presented in \cite{chamarthi2019first}. More recently, we derived discontinuous \cite{wimmer2023fast} and continuous Galerkin \cite{wimmer2024fast} finite element formulations based on a certain mixed formulation, incorporating stabilization techniques used in hyperbolic equations for directional gradients.


For general meshes and implicit time integration, fast linear solvers also pose a major challenge for extreme levels of anisotropy. For high-performance simulations, ideally one would use multilevel methods. However, geometric multigrid methods require some combination of line/plane relaxation and semi-coarsening \cite{Trottenberg_etal_2001,Briggs-etal-2000}, which is expensive, difficult to realize in parallel, and limited to structured grids. Algebraic multigrid (AMG) methods can automatically semi-coarsen in non-grid-aligned directions (e.g., see \cite[Fig. 3]{sivas2021air}), and specialized AMG methods have been developed for anisotropic diffusion, e.g. \cite{schroder2012smoothed,manteuffel2017root,gee2009new,Brandt-etal-2015}, but we have found that these methods still fail on realistic non-grid-aligned anisotropies (e.g., \cite{wimmer2023fast}). In \cite{wimmer2023fast} we develop a specialized mixed discretization and block preconditioner for highly anisotropic diffusion built on AMG for hyperbolic transport operators \cite{manteuffel2018nonsymmetric,manteuffel2019nonsymmetric}, but the solver studies therein are limited to open (i.e. acyclic) field lines, which is not realistic in many practical fusion settings.

In order to reduce the size of the system and construct a computationally efficient solver, we propose a coarse grid approximation built on spectral multiscale basis functions. The concept of using local spectral problems to accurately approximate highly heterogeneous media was proposed in the Generalized Multiscale Finite Element Method (GMsFEM) \cite{efendiev2011multiscale, efendiev2013generalized}. Later, this method was applied to many different applications, such as transport and flow problems, pororoelastity problems in heterogeneous and fractured media, electro-chemical processes in Lithiun-Ion Batteries, geothermal research simulations and seismic wave propagations \cite{chung2018multiscale, vasilyeva2018multiscale, akkutlu2018multiscale, vasilyeva2019upscaling, vasilyeva2021multiscale, vasilyeva2024decoupled}. 
Recently, we adapted the spectral multiscale approach to a general class of discrete problems described using the graph Laplacian \cite{vasilyeva2024generalized, vasilyeva2024adaptive}. The proposed approach is based on spectral characteristics of the normalized local graph Laplacian. 
The central idea of the GMsFEM was associated with a multicontinuum approach and used to construct a nonlocal multicontinuum upscaling (NLMC) \cite{chung2018constraint, chung2018non}. The NLMC is based on constrained energy minimization and provides a good approximation on a coarse grid with sufficient oversampling layers in the basis construction. Moreover, it has been shown that constructing a prolongation operator based on the solution of the local spectral problems gives a basis for constructing highly efficient preconditioners in the spectral algebraic multigrid method \cite{chartier2003spectral, falgout2005two, brezina2011smoothed, galvis2010domain, efendiev2012multiscale}. Most traditional approaches in multiscale model order reduction address highly varied properties in a classic manner related to heterogeneous properties, but do not explore high-contrast variations caused by physical fields that span the entire domain. Some work in element-based AMG has been successfully applied to anisotropic diffusion, e.g. \cite{chartier2007spectral}, but the problems considered have sufficiently mild anisotropy that existing AMG solvers (i.e., without spectral coarse modes), e.g. \cite{manteuffel2017root,schroder2012smoothed,gee2009new,Brandt-etal-2015}, can still be effective solvers. In this work, we address the construct associated with the two heat flow directions: parallel and perpendicular. As mentioned earlier, in real applications, the ratio between parallel and perpendicular heat conductivities can reach $10^{10}$ or higher, and the main challenge is associated with the parallel direction of heat flow following the magnetic field, particularly for non-grid-aligned meshes. 

In this work, we show that the local generalized spectral problems applied to the anisotropic diffusion equation yield eigenfunctions that align with the magnetic field lines and give a very good approximation of the problem on a coarse grid. Thus, even without enforcing mesh alignment with a given magnetic field, we are able to integrate in a basis that effectively \emph{does} respect the magnetic field structure.
More specifically, we start with a finite element approximation on a sufficiently fine grid with second-order polynomial basis functions, and use this as a reference solution to estimate the accuracy of the proposed coarse grid approximation. We then construct a highly accurate approximation by introducing a local generalized eigenvalue problem and use eigenvectors associated with the smallest eigenvalues to construct a multiscale space aligned with a magnetic field-directed heat flux. We show approximation properties of the constructed multiscale space, and then incorporate the coarse multiscale space into a two-grid preconditioner. Multiscale approximation results are used to prove convergence of the two-grid method, independent of anisotropy. We present numerical accuracy investigations for three test problems with different magnetic fields and varying contrasts of parallel and perpendicular heat conductance. We show that with sufficiently many multiscale basis functions, we can obtain highly accurate results with a large system size reduction, and in the context of a two-grid preconditioner, achieve convergence in $\mathcal{O}(1)$ iterations, independent of anisotropy. 

The paper is structured as follows. In Section \ref{sec:problem}, we outline the problem formulation for anisotropic diffusion, presenting weak formulations and approximations on a fine grid. Section \ref{sec:ms} presents a coarse grid approximation using a spectral multiscale space and discusses approximation properties of the proposed multiscale method. A two-grid method built around the spectral multiscale coarse space is discussed in Section \ref{sec:tg}, including a proof of two-grid convergence. Numerical results are presented in Section \ref{sec:results} to illustrate the robustness and accuracy of the proposed method for three test cases with different magnetic field setups. We show that the method provides an accurate approximation given a sufficient number of multiscale basis functions and coarse grid size. The paper ends with a conclusion.

\section{Problem formulation}\label{sec:problem}
 
We consider anisotropic diffusion in a domain $\Omega$ for temperature $T = T(t, x)$:
\begin{equation}
\label{eq:1a}
T_t - \nabla \cdot \mathbf{q} = f, \quad
x \in \Omega, \quad t > 0,  
\end{equation}
with 
\begin{equation}
\label{eq:1b}
\mathbf{q}  = k_{\parallel} \nabla_{\parallel} T + k_{\perp} \nabla_{\perp} T, \quad 
\nabla_{\parallel} (\cdot) = ( \mathbf{b} \cdot \nabla (\cdot)) \mathbf{b}, \quad 
\nabla_{\perp} = \nabla - \nabla_{\parallel},
\end{equation}
with given initial conditions $T(0, x) = T_0$ for $x \in \Omega$ and Dirichlet boundary conditions
\[
T(t, x) = g, \quad x \in \partial \Omega, \quad t > 0.
\]
Here $\mathbf{b} = \mathbf{B}/|\mathbf{B}|$, where $\mathbf{B}$ is the given magnetic field, and $k_{\parallel}$ and $k_{\perp}$ are the heat conductivity parallel and perpendicular to normalized magnetic field lines $\mathbf{b}$, respectively. 

Substituting \eqref{eq:1b} to \eqref{eq:1a}, we obtain the following formulation: 
\begin{equation}
\label{eq:1}
T_t - \nabla \cdot (k_{\perp} \nabla T)
- \nabla \cdot ( k_{\Delta} \ \mathbf{b} (\mathbf{b} \cdot \nabla T) )
 = f, 
\end{equation}
with $k_{\Delta} = k_{\parallel} - k_{\perp}$.

Let
$V = \{v \in H^1(\Omega): v = g \ \text{on} \ \partial \Omega \}$ and  
$\hat{V}  = \{v \in H^1(\Omega): v = 0 \ \text{on} \ \partial \Omega \}$. 
Then, we can write the following variational formulation: find $T \in V$ such that
\begin{equation}
\label{eq:2}
m (T_t, v)   + a(T, v) = l(v), \quad \forall v \in \hat{V},
\end{equation}
with 
\begin{equation}
\label{eq:2a}
\begin{split}
&a(T, v) = a_{iso}(T, v) +  a_{aniso}(T, v), \quad 
m(T, v)  = \int_{\Omega} T \ v\ dx, \quad 
l(f) = \int_{\Omega} f \ v\ dx, \\
&
a_{iso}(T, v)  = \int_{\Omega} k_{\perp} \nabla T \cdot \nabla v \ dx, \quad 
a_{aniso}(T, v)  = \int_{\Omega} k_{\Delta}  (\mathbf{b} \cdot \nabla T)  (\mathbf{b} \cdot \nabla v) \ dx.
\end{split}
\end{equation}

Let $\mathcal{T}_h$ be a finite element partition of the domain into elements $K_i$ with mesh size $h$
\[
\mathcal{T}_h = \cup_{i=1}^{N_h^{cell}} K_i,
\] 
where $N_h^{cell}$ is the number of cells.  We call $\mathcal{T}_h$ a ``fine grid,'' where we assume the fine grid has sufficient resolution of the magnetic field $\mathbf{b}$ and temperature field $T$ that it can provide an accurate approximate solution using standard polynomial basis functions. 
For approximation space on the fine grid $\mathcal{T}_h$, we use a continuous Galerkin (CG) formulation. 
The variational problem reads: find $T \in V_h \subset V$ such that
\begin{equation}
\label{eq:3}
m \left((T_h)_t, v \right)  + a(T_h, v) = l(v), \quad \forall v \in \hat{V}_h \subset \hat{V},
\end{equation}
where 
\[
T_h = \sum_{j=1}^{N_h} T_j \phi_j(x), \quad 
V_h = \text{span} \{\phi_j, \ j=1, \ldots N_h \}.
\]

Let $\tau$ be the time step size and $T_h^n$ be the solution at time $t^n$. 
Then, using a backward Euler time approximation, we obtain the following discrete problem: find $T^n_h \in V_h$ such that
\begin{equation}
\label{eq:5}
\frac{1}{\tau} m \left( T^n_h , v \right) 
+ a (T^n_h, v)  = \frac{1}{\tau} m \left(T^{n-1}_h, v \right)  + l(v), \quad \forall v \in \hat{V}.
\end{equation}
Note, this form generalizes to implicit multistage or multistep methods as well, with a modified linear forcing function. Altogether, we arrive at a discrete problem for each implicit time step \eqref{eq:5} in the following matrix form
\begin{equation}
\label{eq:6}
\left(\frac{1}{\tau} M_h + A_h \right) T^n_h  = 
\frac{1}{\tau} M_h T^{n-1}_h + F_h,
\end{equation}
where 
$M_h = \{m_{ij} = m(\phi_i, \phi_j)\}$, 
$A_h = \{a_{ij} = a(\phi_i, \phi_j)\}$ and 
$F_h = \{f_j = l(\phi_j)\}$.

\section{Multiscale space approximation}\label{sec:ms}

\begin{figure}[h!]
\centering
\includegraphics[width=0.48\linewidth]{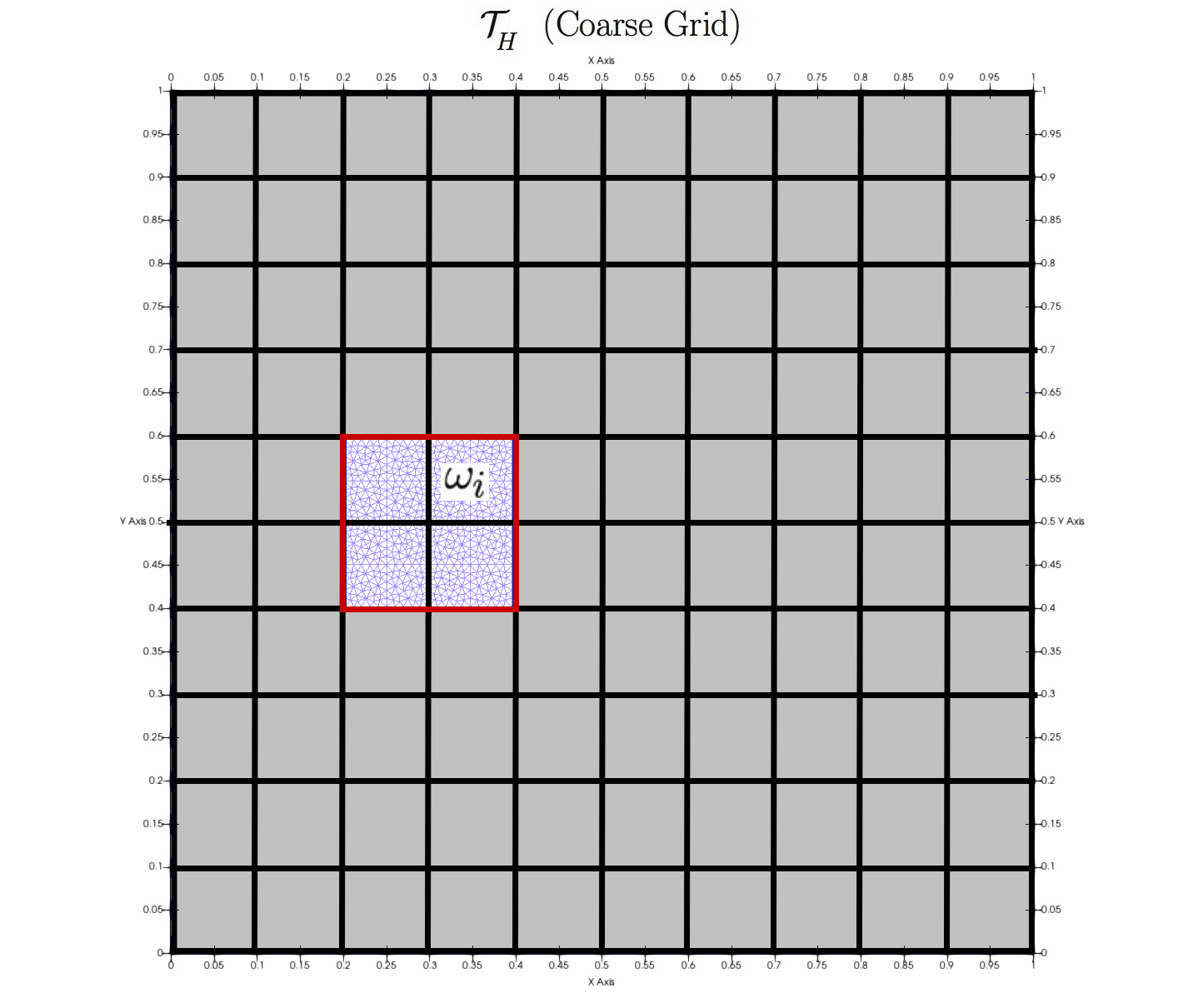}
\includegraphics[width=0.48\linewidth]{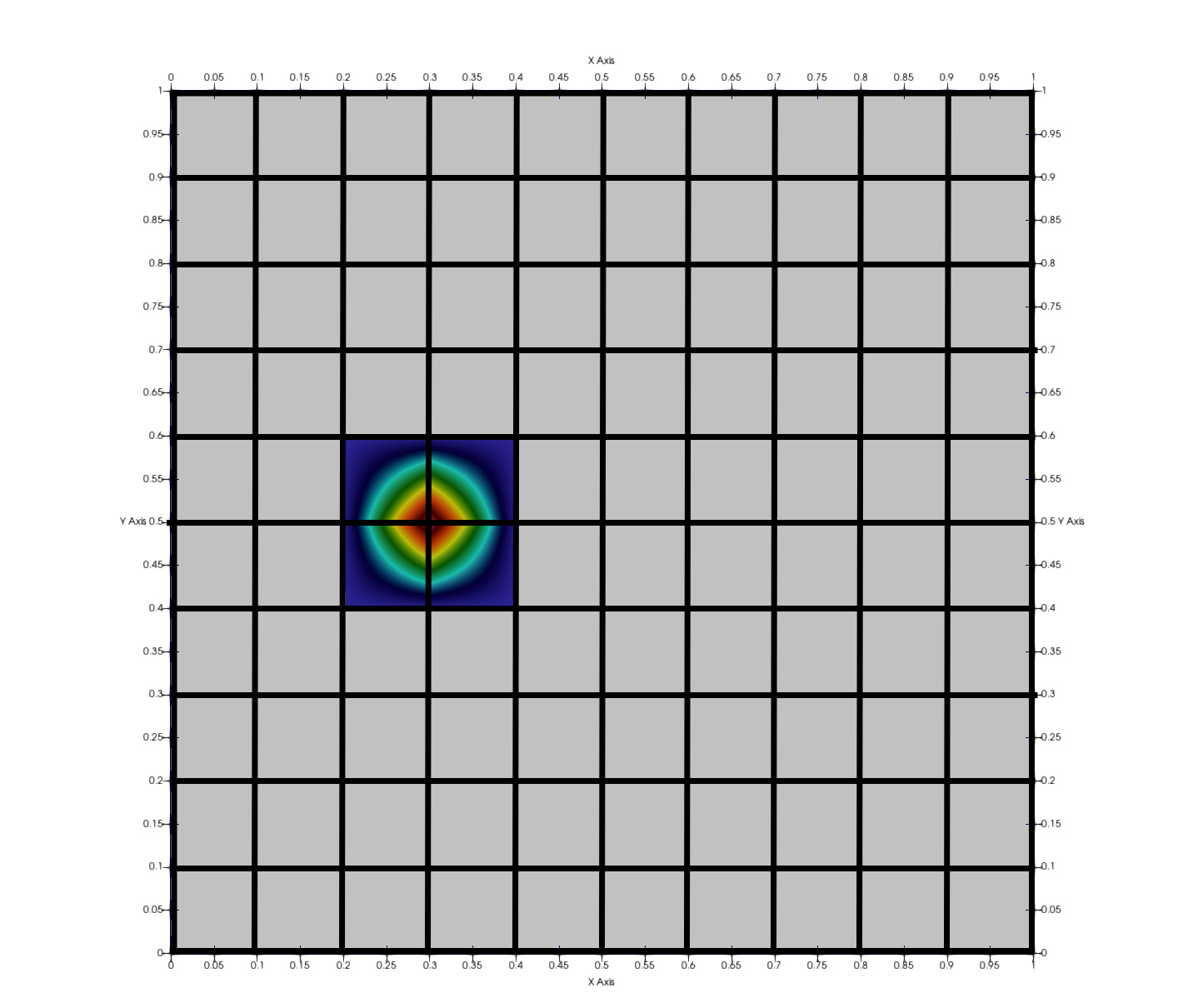}
\caption{Illustration of the $10 \times 10$ coarse grid $\mathcal{T}_H$ with local domain $\omega_i$ and linear partition of unity function $\chi_i$}
\label{meshc}
\end{figure}

Let $\mathcal{T}_H$ be a finite element partition of the domain and mesh into coarse elements $\{K_i\}$ with mesh size $H >> h$:
\[
\mathcal{T}_H = \cup_{i=1}^{N_H^{cell}} K_i,
\] 
where $N_H^{cell}$ is the number of coarse cells (left plot in Figure \ref{meshc}). We then consider a CG formulation for a coarse grid approximation as well, and construct nodal multiscale basis functions for the coarse variational problem. 

For nodal basis functions, let $\omega_i$ be the subdomain defined as the collection of coarse cells containing coarse grid node $x_i$. Let $\mathcal{T}^{\omega_i}_h$ be the fine-scale partitioning of the local domain $\omega_i$. 
We suppose that $\mathcal{T}^{\omega_i}_h$ is conforming with global finite element partitioning $\mathcal{T}_h$, so we are able to build a local-to-global map for all local degrees of freedom $DOF_h^{\omega_i} = N_h^{\omega_i}$ to global $DOF_h = N_h$. 
In addition, we define a partition of unity function $\chi_i$ defined in $\omega_i$. This approach will use a regular bilinear partition of unity functions (right plot in Figure \ref{meshc}). 

For the construction of a basis that provides an accurate approximation on the coarse grid $\mathcal{T}_H$, we will construct spectral multiscale basis functions motivated by the Generalized Multiscale Finite Element Method (GMsFEM) \cite{efendiev2011multiscale, efendiev2013generalized, chung2016generalized}. In the presented approach, we design an eigenvalue problem based on the normalized anisotropic Laplacian operator recently presented for a general class of problems on graph \cite{vasilyeva2024generalized}. We introduce a local spectral problem in each local domain $\omega_i$ and show that the eigenvectors naturally follow the underlying magnetic field lines, and can provide superior approximation for the considered problem compared with typical piecewise polynomials.

Assuming a fixed in time magnetic field, the algorithm for the multiscale method contains offline and online stages:
\begin{itemize}
\item \textit{Offline calculations. }
\begin{itemize}
\item Define coarse grid $\mathcal{T}_H$ and generate local domains $\omega_i$.
\item Solve local spectral problems to construct a set of multiscale basis functions $\{\psi^{\omega_i}_j\}$ in each local domain $\omega_i$ independently.
\item Map the local degrees of freedom to global and form an interpolation operator $P$ using a given number of $J$ local multiscale basis functions over each subdomain, where
\[
P = \left[ 
\psi^{\omega_1}_1, \ldots, \psi^{\omega_1}_{J},
\ldots
\psi^{\omega_{N_H^{vert}}}_1, \ldots, \psi^{\omega_{N_H^{vert}}}_{J}
\right].
\]
\end{itemize}
\item \textit{Online calculations. }
\begin{itemize}
\item Project matrices and vectors to the coarse grid using the precomputed interpolation matrix $P$ and restriction matrix $R = P^T$.
\item Solve the coarse-scale system.
\item Interpolate the coarse-scale solution to the fine grid resolution.
\end{itemize}
\end{itemize}
Next, we describe the construction of the spectral multiscale basis functions and discuss approximation properties of the multiscale space.

\subsection{Spectral multiscale basis functions}

Within the fine grid resolution of the local domain $\omega_i$, we solve the following generalized eigenvalue problem in matrix form
\begin{equation} 
\label{eq:sp}
A^{\omega_i}_h \phi^{\omega_i} = 
\lambda D^{\omega_i}_h \phi^{\omega_i}, 
\end{equation}
where
\[
A^{\omega_i}_h = \{a^{\omega_i}_{ij} = a^{\omega_i}(\phi_i, \phi_j)\}, \quad 
a^{\omega_i}(u, v) = a^{\omega_i}_{iso}(u, v) +  a^{\omega_i}_{aniso}(u, v), 
\]\[
a^{\omega_i}_{iso}(u, v)  = \int_{\omega_i} k_{\perp} \nabla u \cdot \nabla v \ dx, \quad 
a^{\omega_i}_{aniso}(u, v)  = \int_{\omega_i} k_{\Delta}  (b \cdot \nabla u)  (b \cdot \nabla v) \ dx,
\]
and $D_h^{\omega_i} = \text{diag}(d^{\omega_i}_1, \ldots,d^{\omega_i}_{N^{\omega_i}_h})$, where $d^{\omega_i}_i = a^{\omega_i}_{ii}$ are the diagonal elements of the matrix $A^{\omega_i}_h$ \cite{vasilyeva2024generalized, chartier2003spectral, falgout2005two, brezina2011smoothed}. Here $A_h^{\omega_i}$ and $D_h^{\omega_i}$ are assumed to be SSPD by nature of the underlying problem and discretization, implying the generalized eigenvalues are real and nonnegative, and the eigenvectors form an orthogonal basis. Moreover, the first eigenvalue is zero with a constant corresponding eigenvector, as only natural boundary conditions are imposed on the local bilinear form.

\begin{figure}[h!]
\centering
\begin{subfigure}[b]{0.32\textwidth}
\centering
\includegraphics[width=1\textwidth]{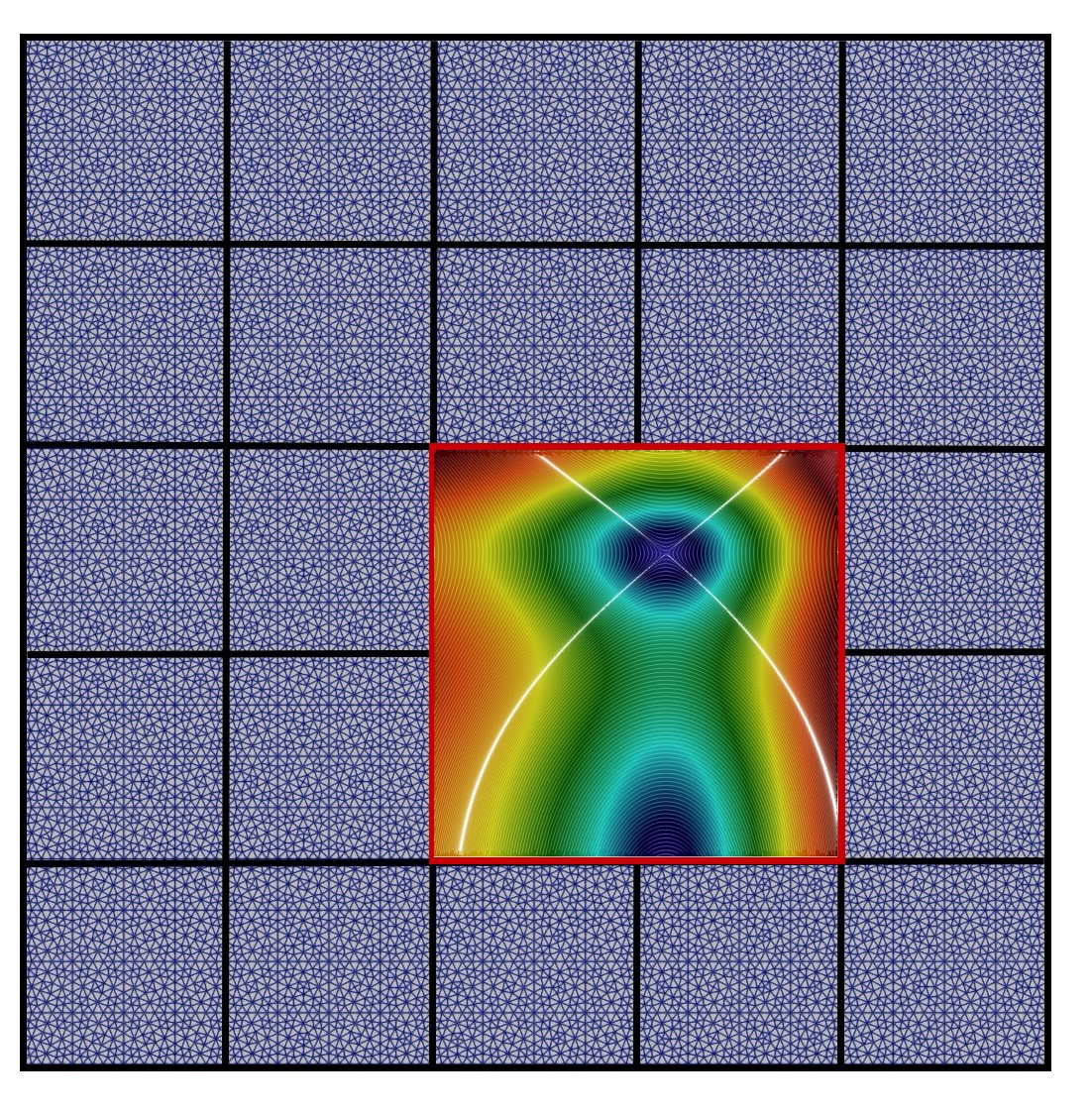} 
\caption{$5 \times 5$ coarse grid and local domain $\omega_{15}$ with magnetic field}
\end{subfigure}
\begin{subfigure}[b]{0.67\textwidth}
\centering
\begin{subfigure}[b]{1\textwidth}
\centering
\includegraphics[width=0.15 \textwidth]{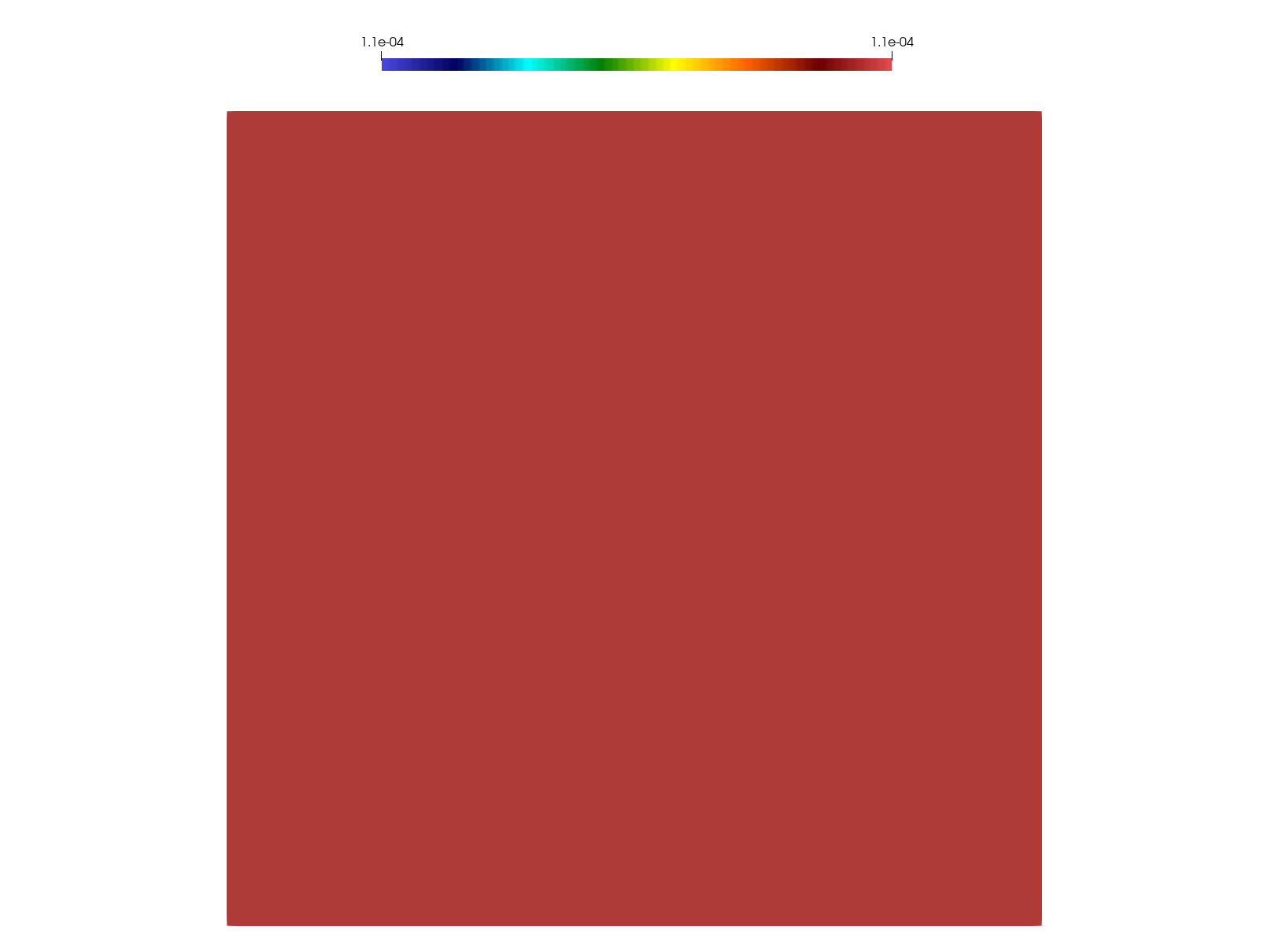} 
\includegraphics[width=0.15 \textwidth]{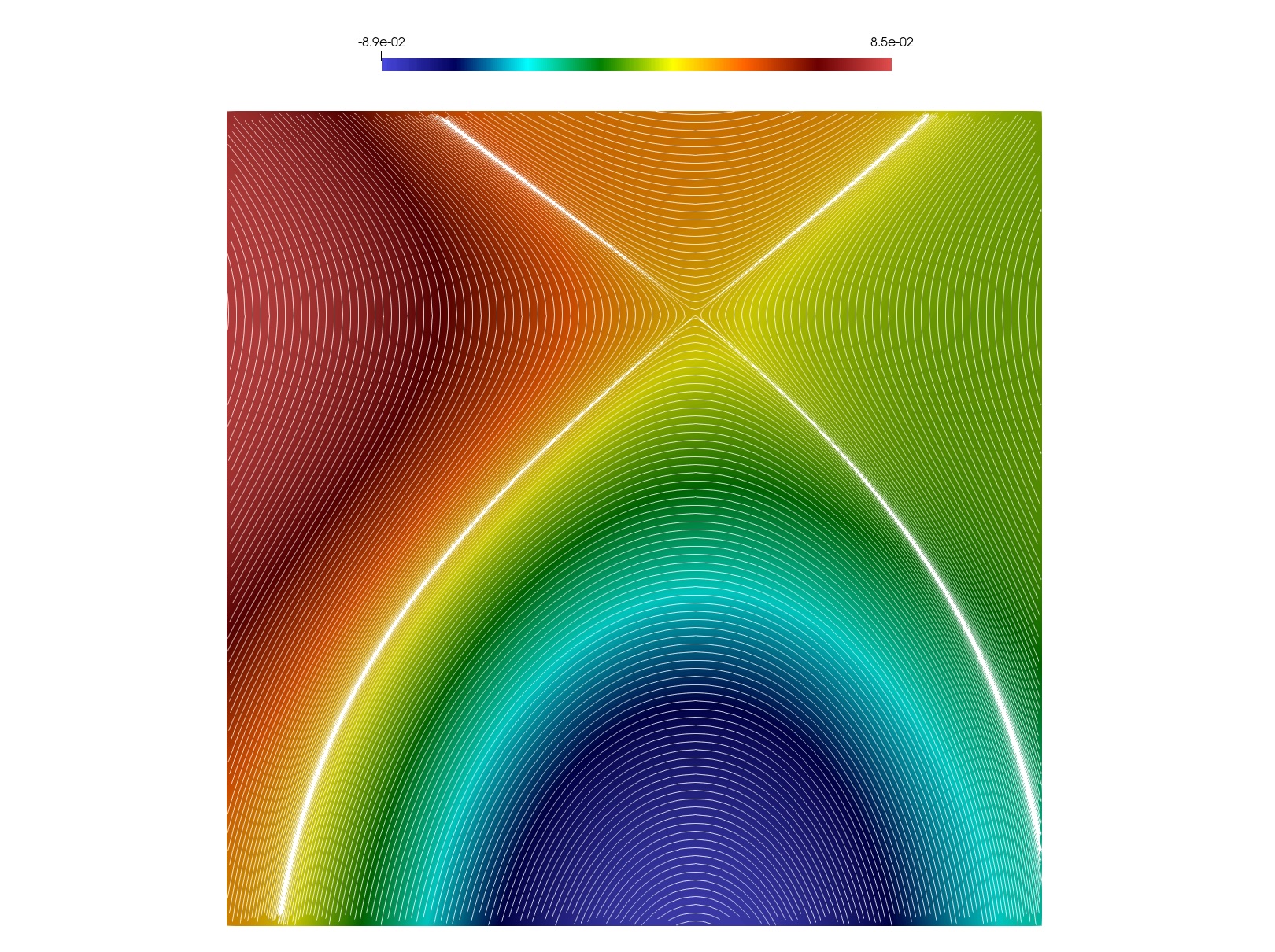} 
\includegraphics[width=0.15 \textwidth]{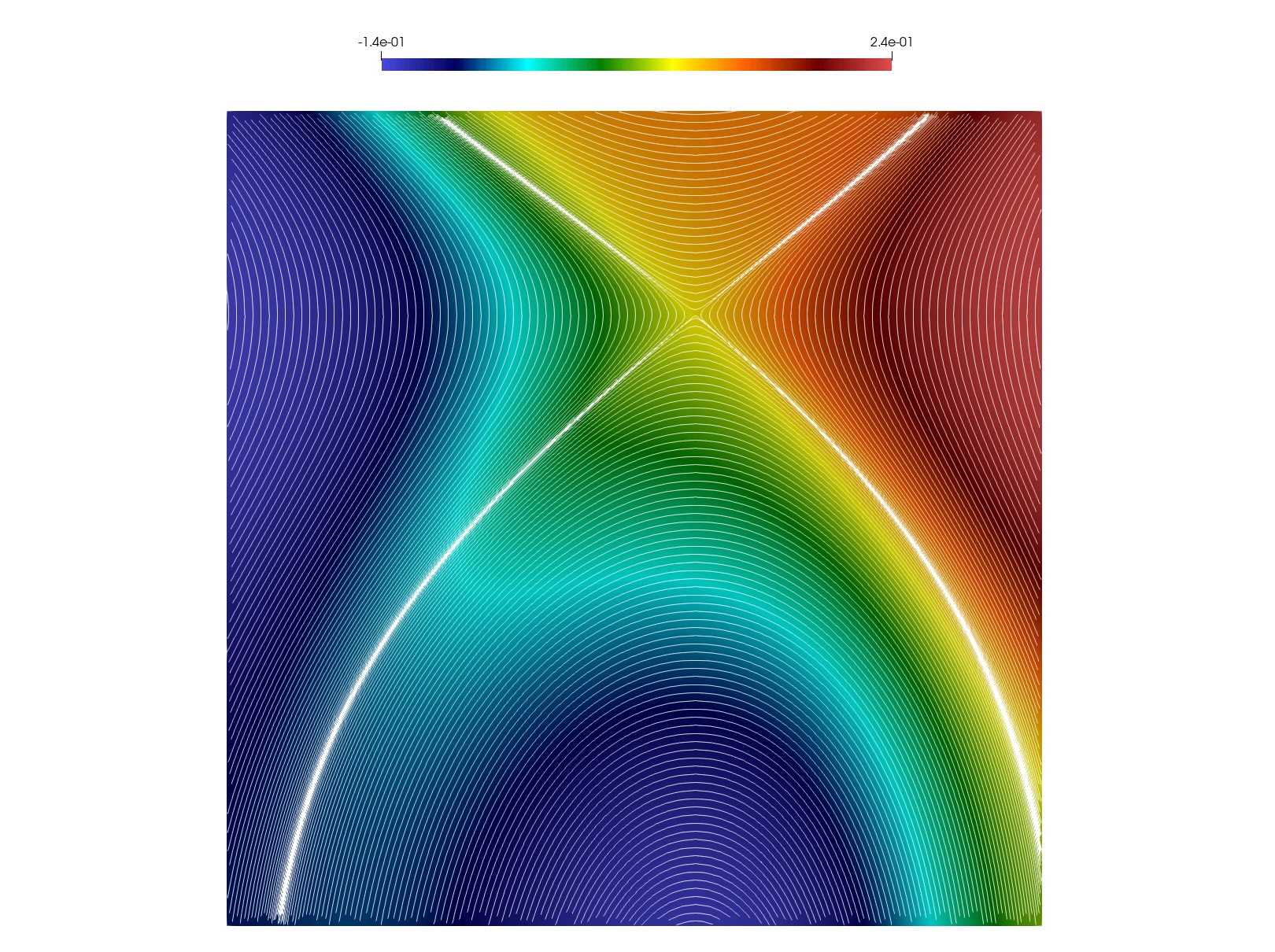} 
\includegraphics[width=0.15 \textwidth]{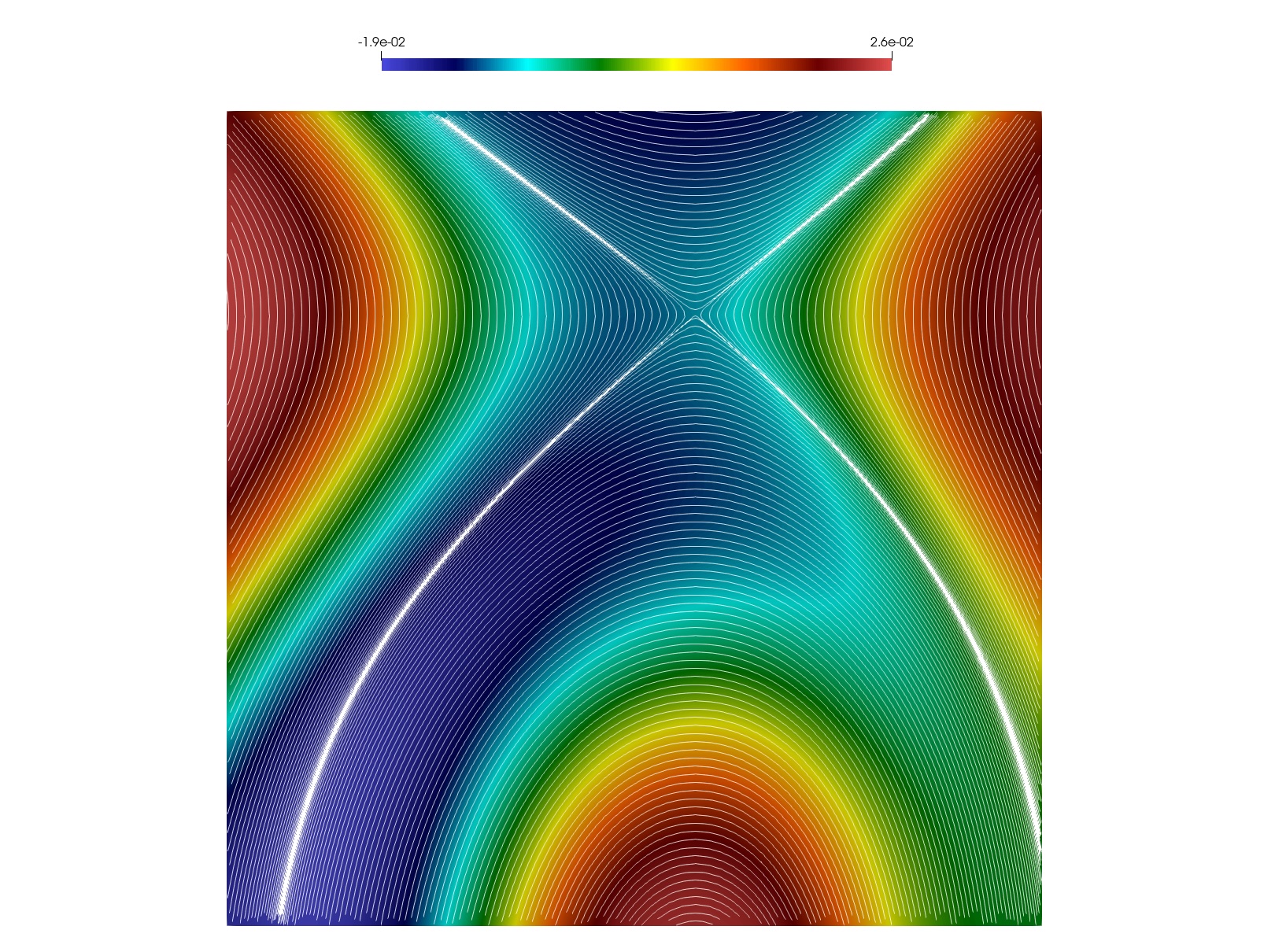} 
\includegraphics[width=0.15 \textwidth]{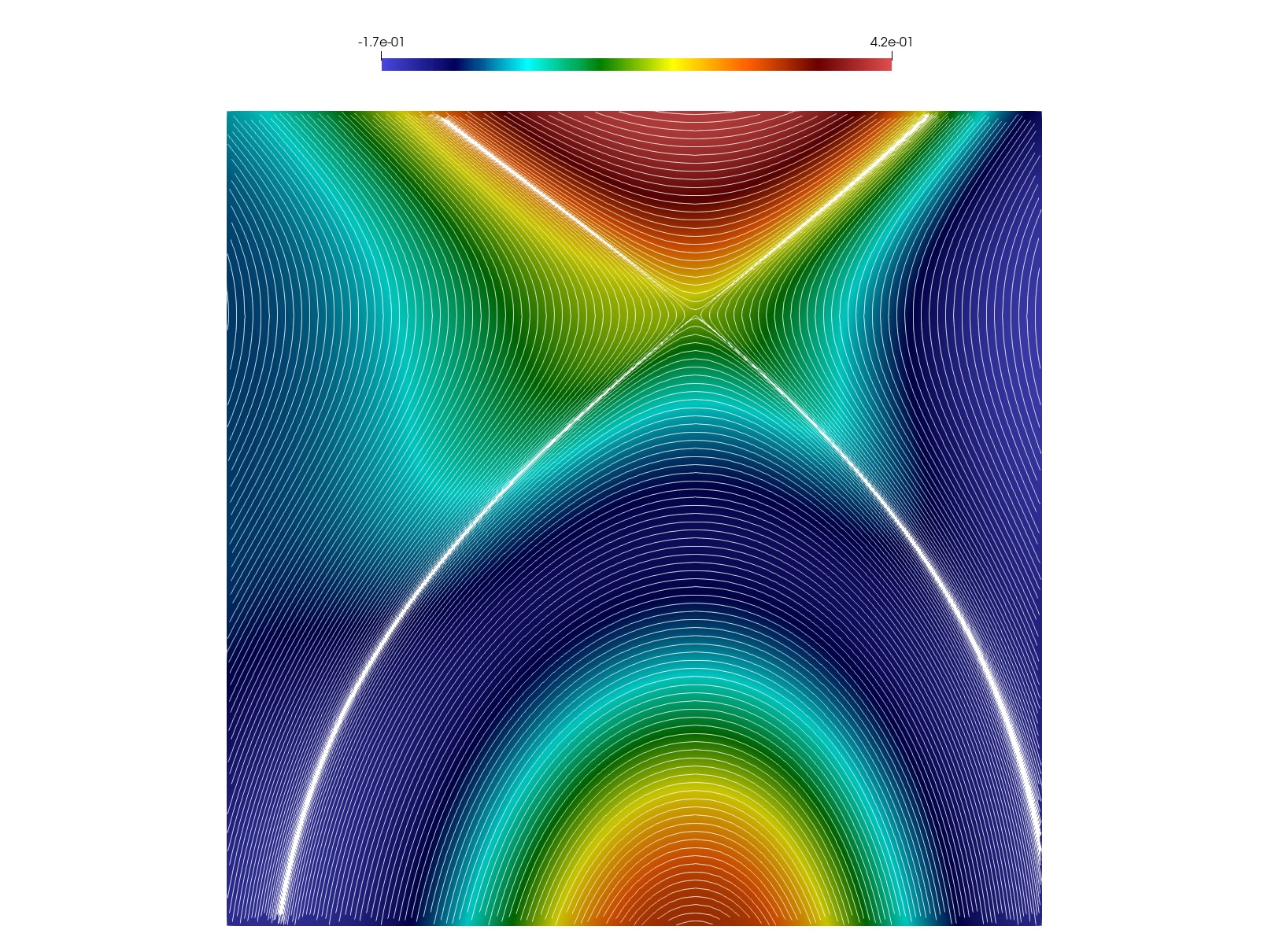} 
\includegraphics[width=0.15 \textwidth]{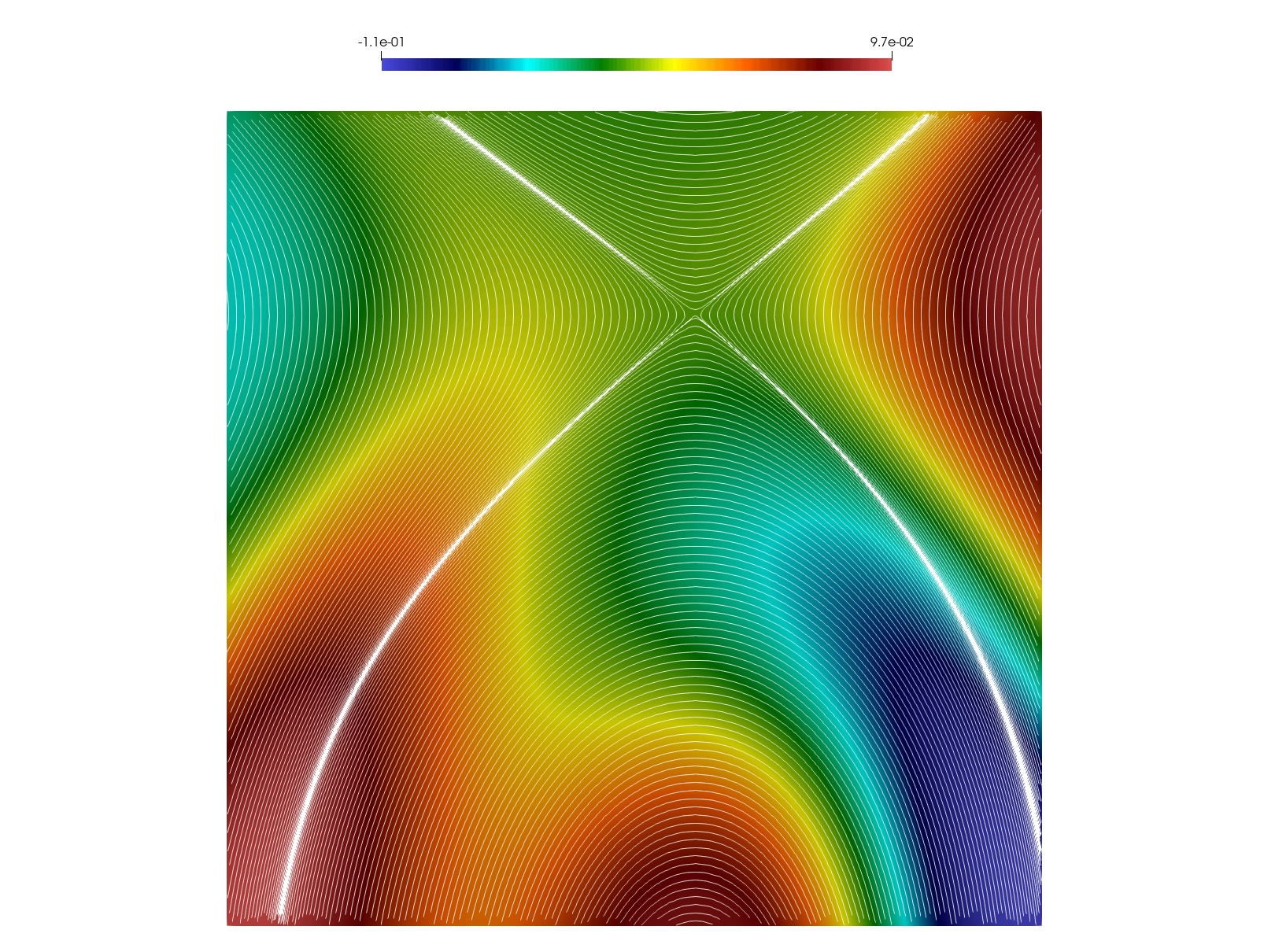} 
\caption{$k_{\parallel}/k_{\perp} = 10^1$}
\end{subfigure}
\begin{subfigure}[b]{1\textwidth}
\centering
\includegraphics[width=0.15 \textwidth]{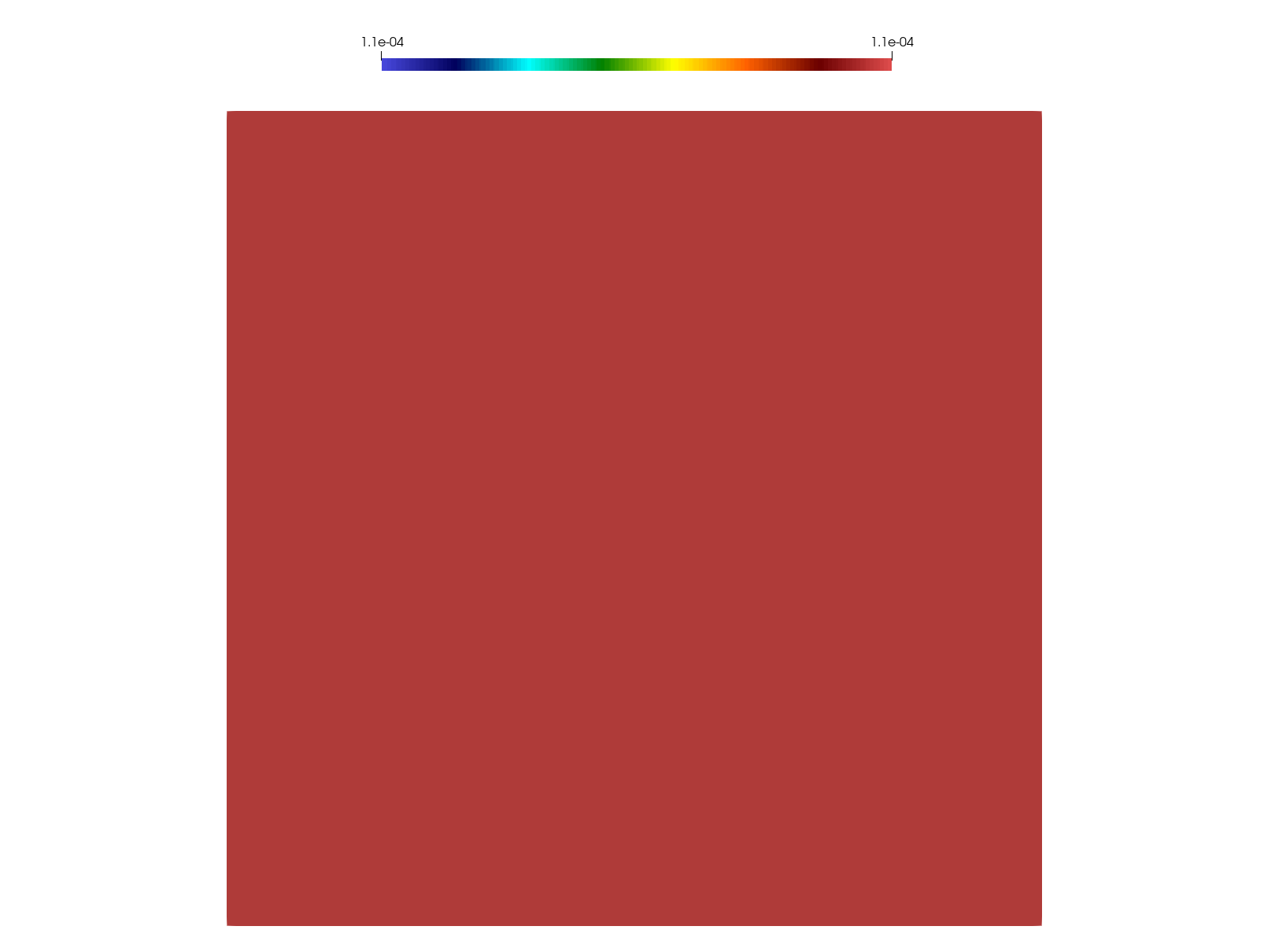} 
\includegraphics[width=0.15 \textwidth]{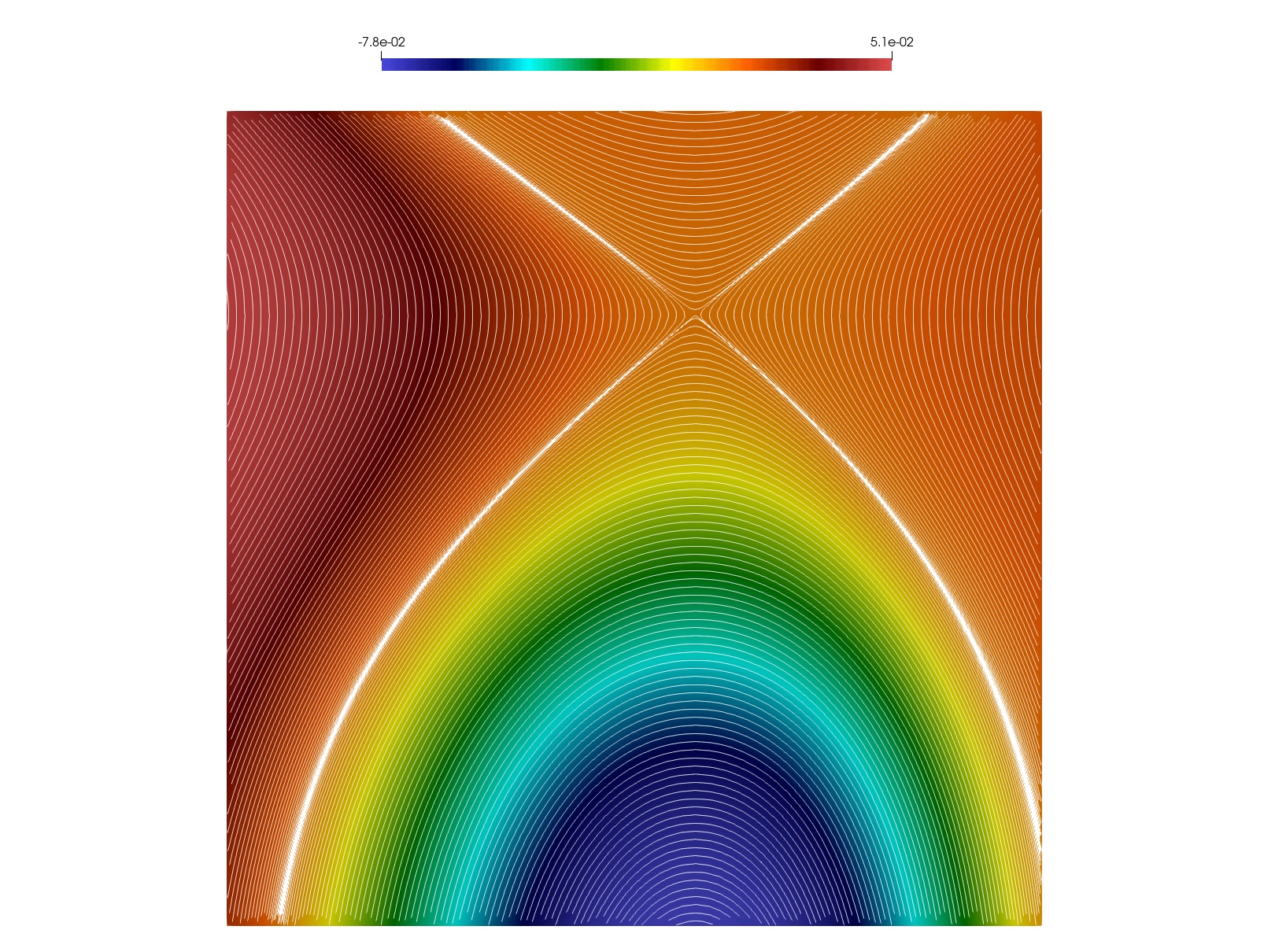} 
\includegraphics[width=0.15 \textwidth]{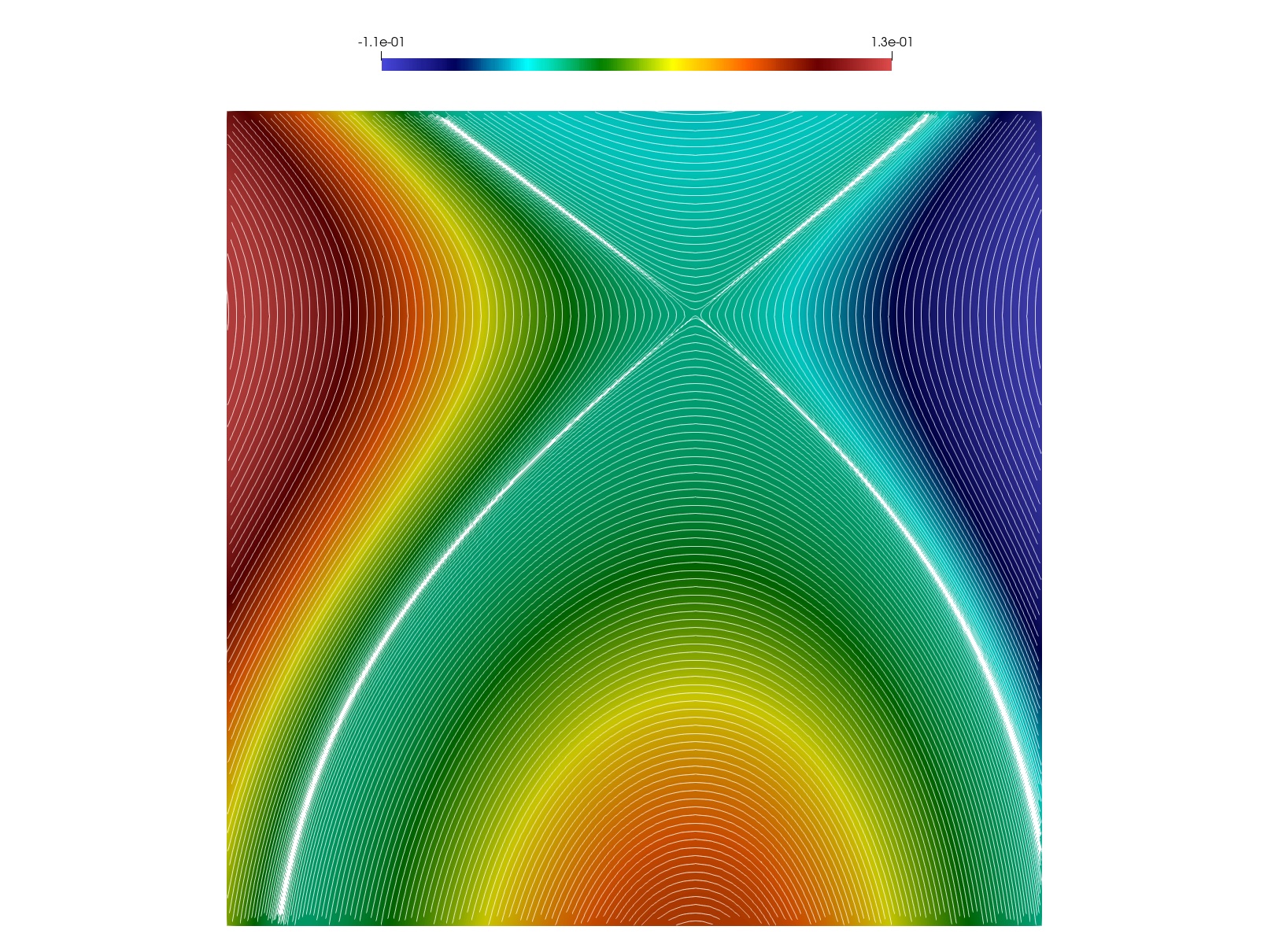} 
\includegraphics[width=0.15 \textwidth]{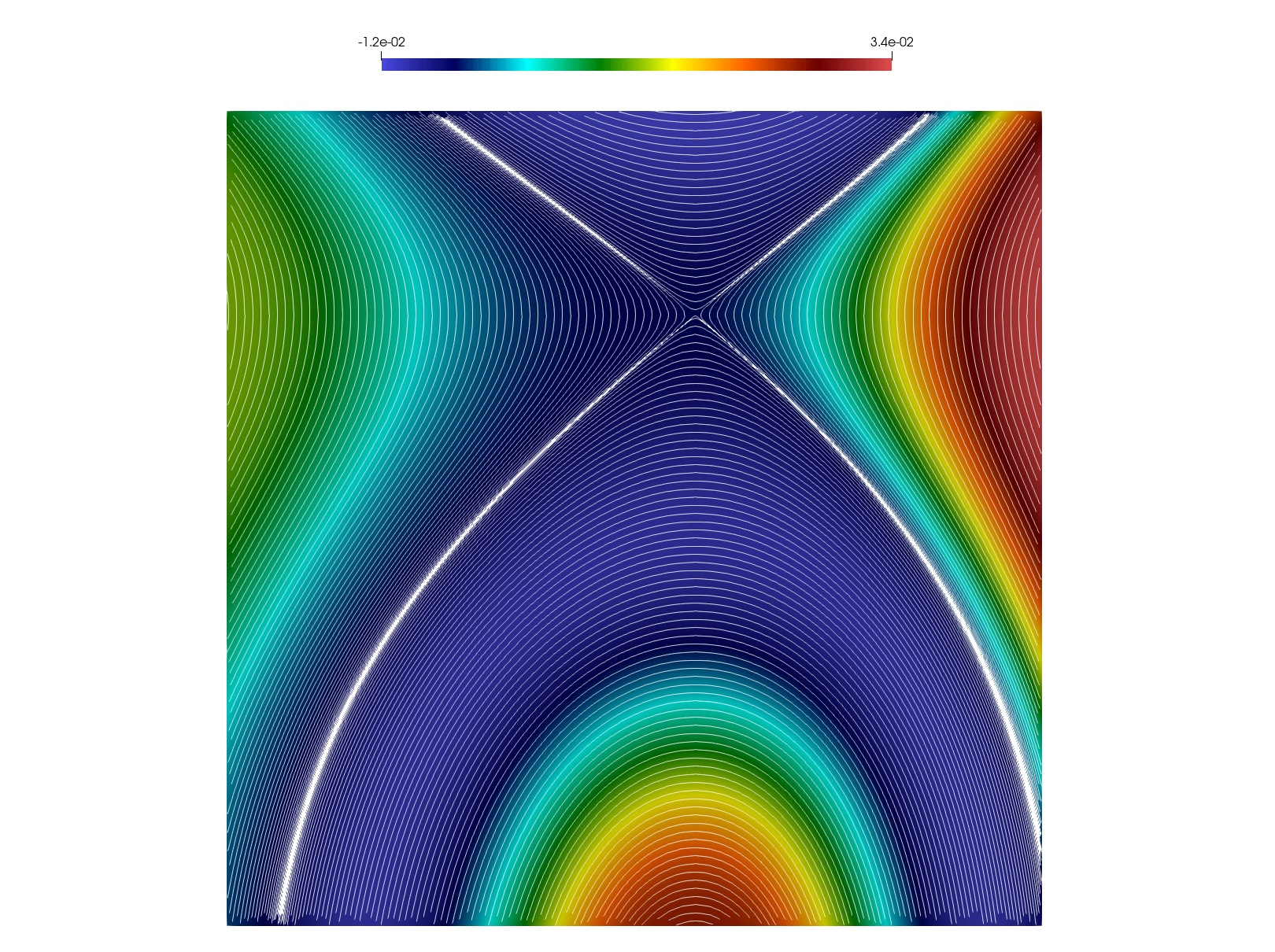} 
\includegraphics[width=0.15 \textwidth]{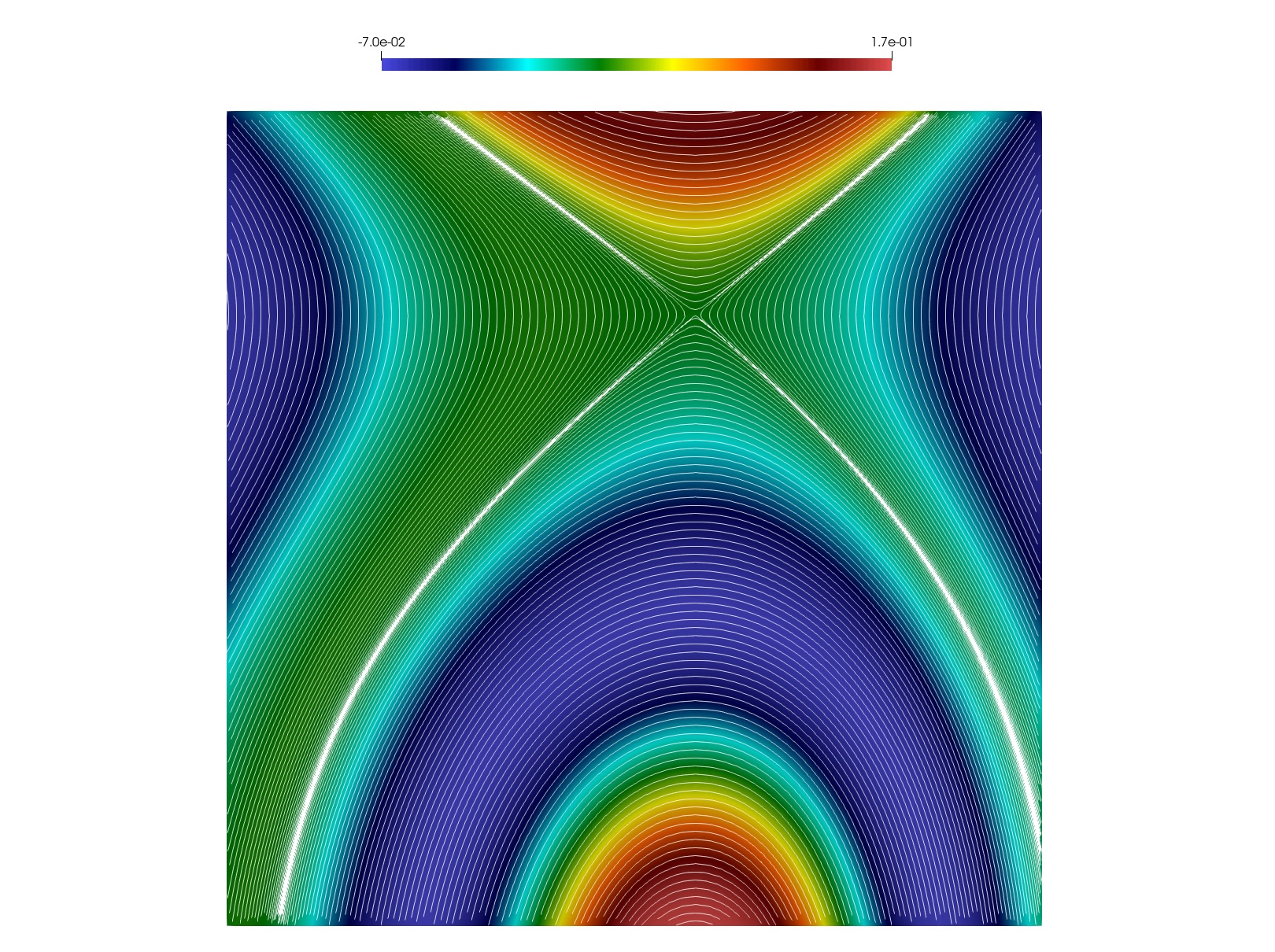} 
\includegraphics[width=0.15 \textwidth]{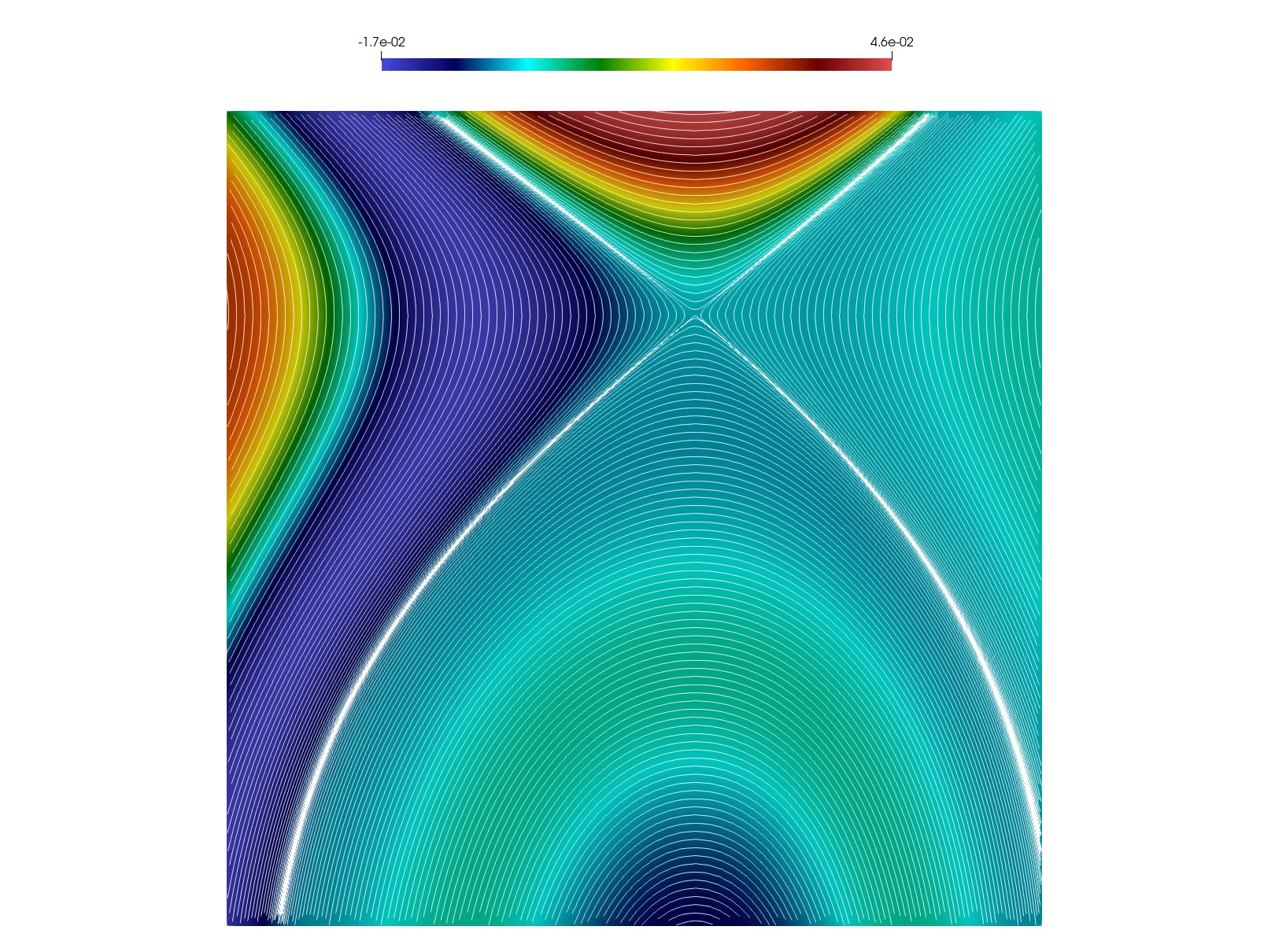} 
\caption{$k_{\parallel}/k_{\perp} = 10^3$}
\end{subfigure}
\begin{subfigure}[b]{1\textwidth}
\centering
\includegraphics[width=0.15 \textwidth]{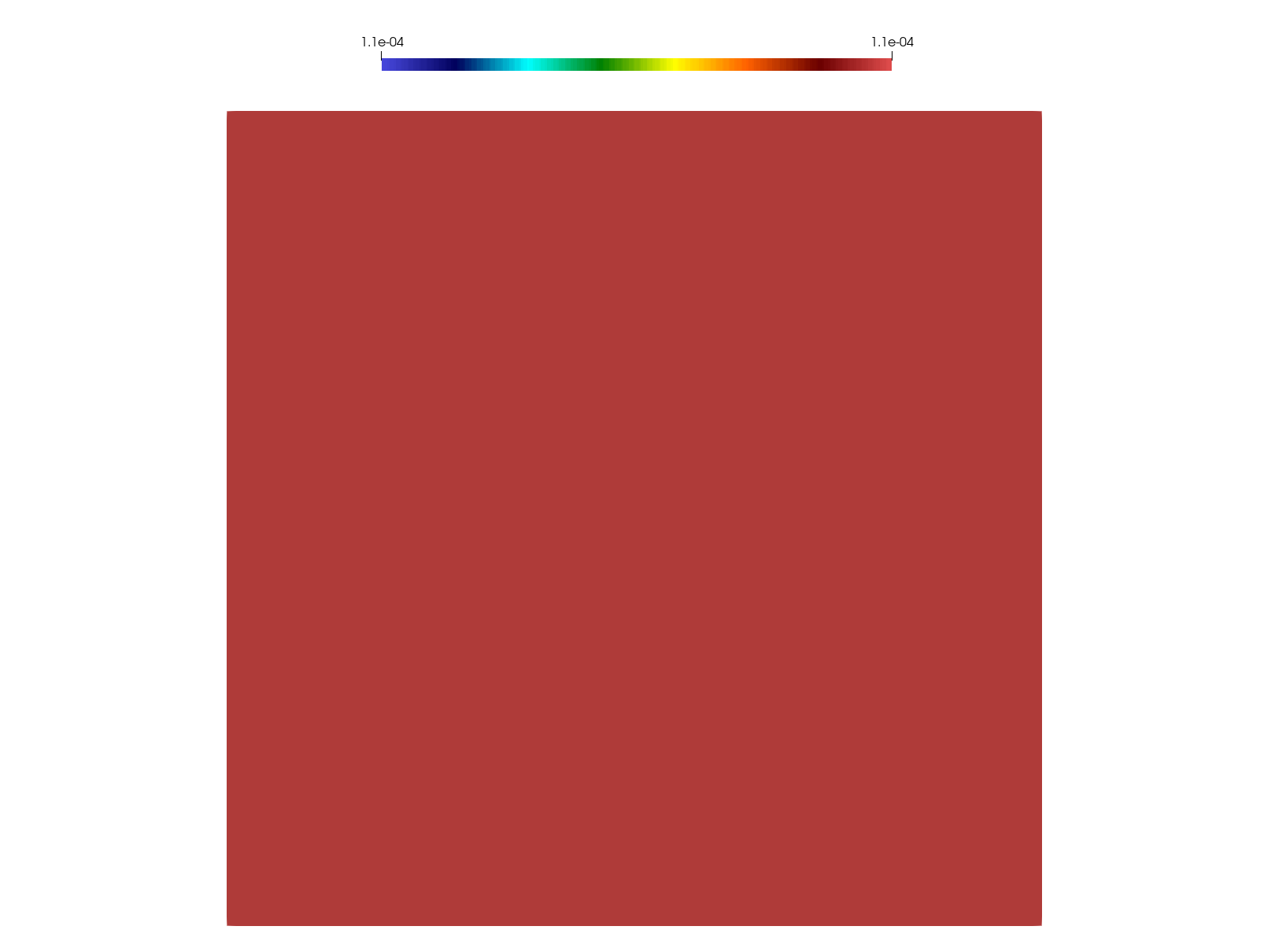} 
\includegraphics[width=0.15 \textwidth]{b/u2} 
\includegraphics[width=0.15 \textwidth]{b/u3} 
\includegraphics[width=0.15 \textwidth]{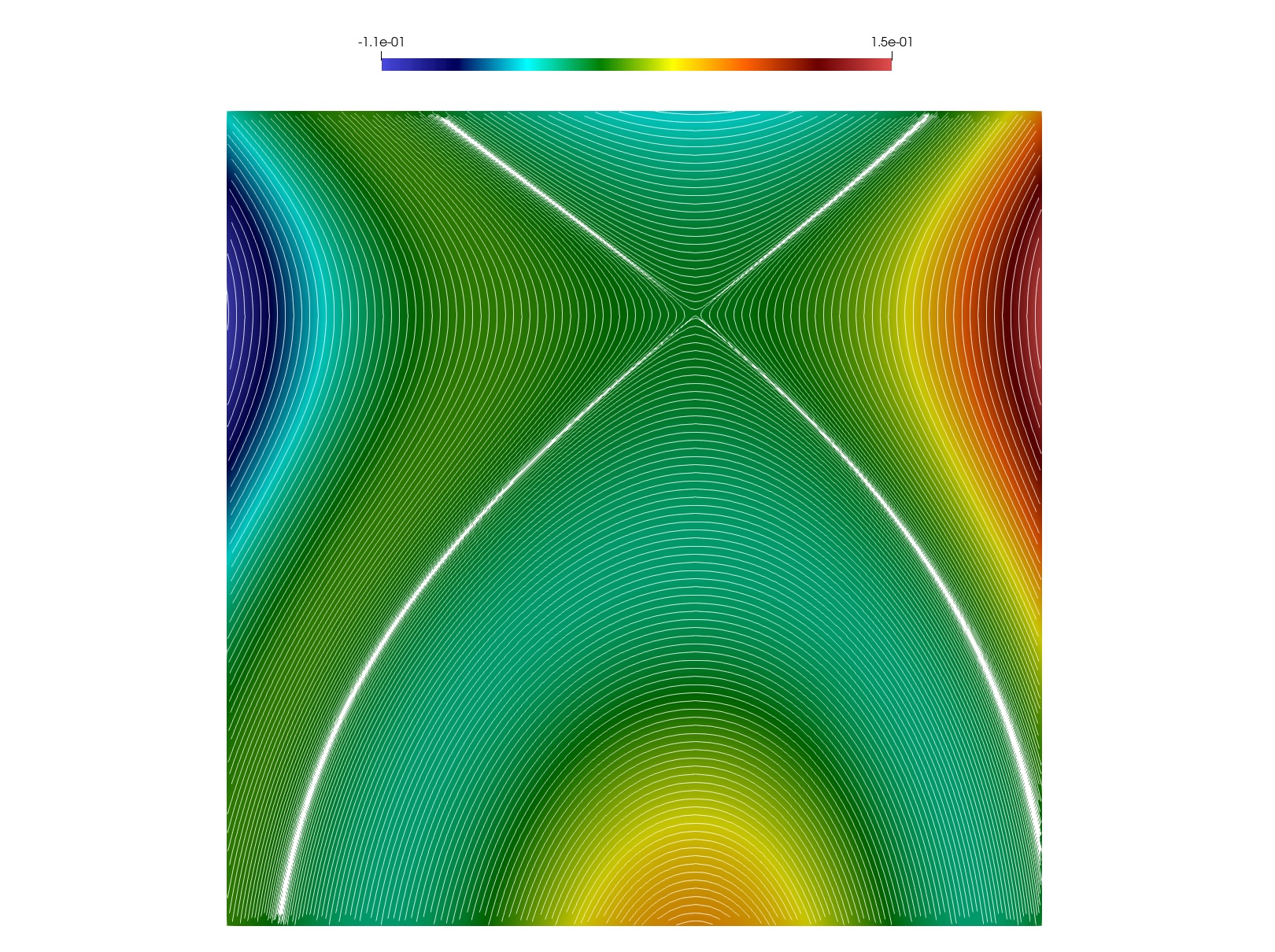} 
\includegraphics[width=0.15 \textwidth]{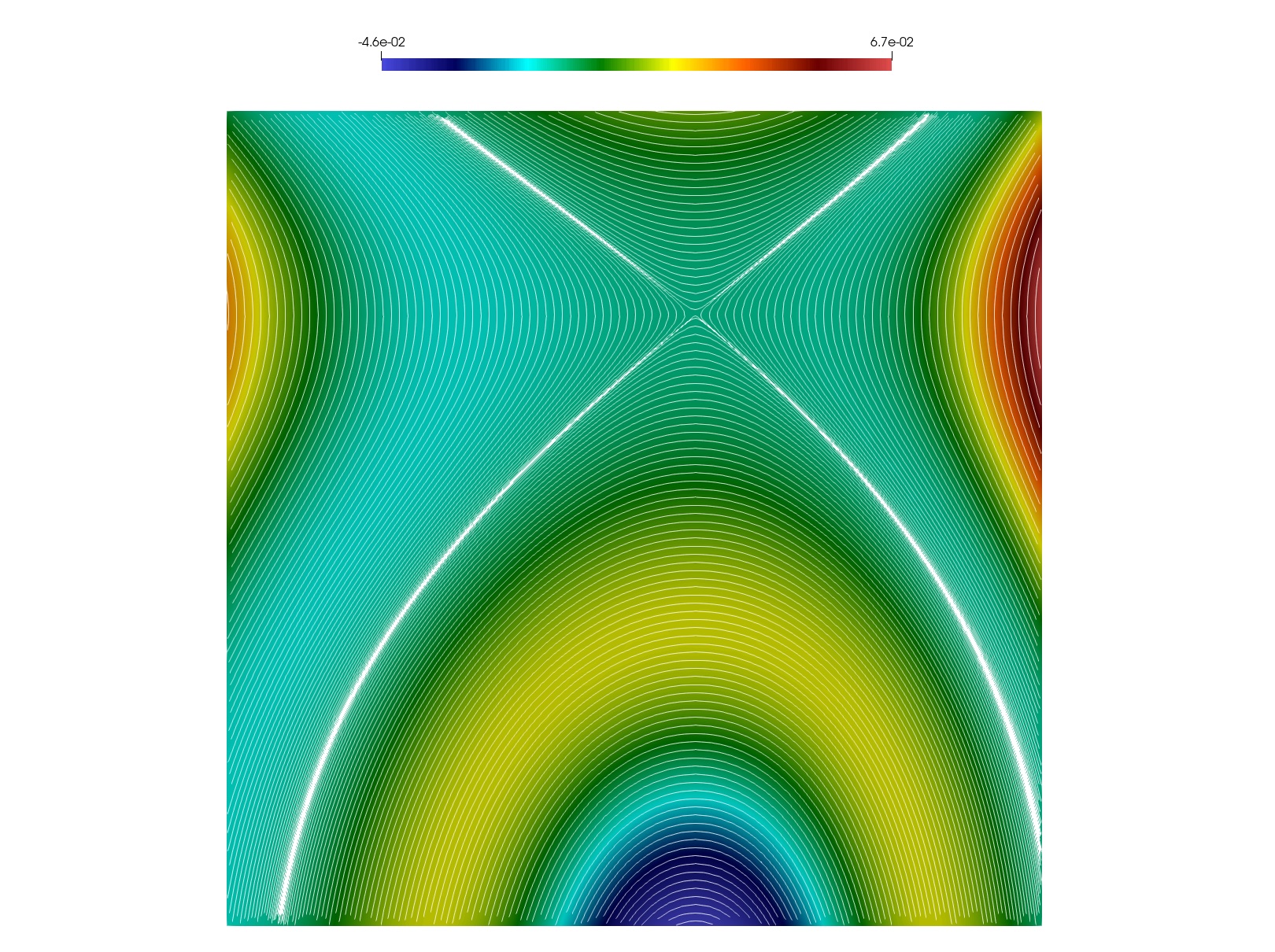} 
\includegraphics[width=0.15 \textwidth]{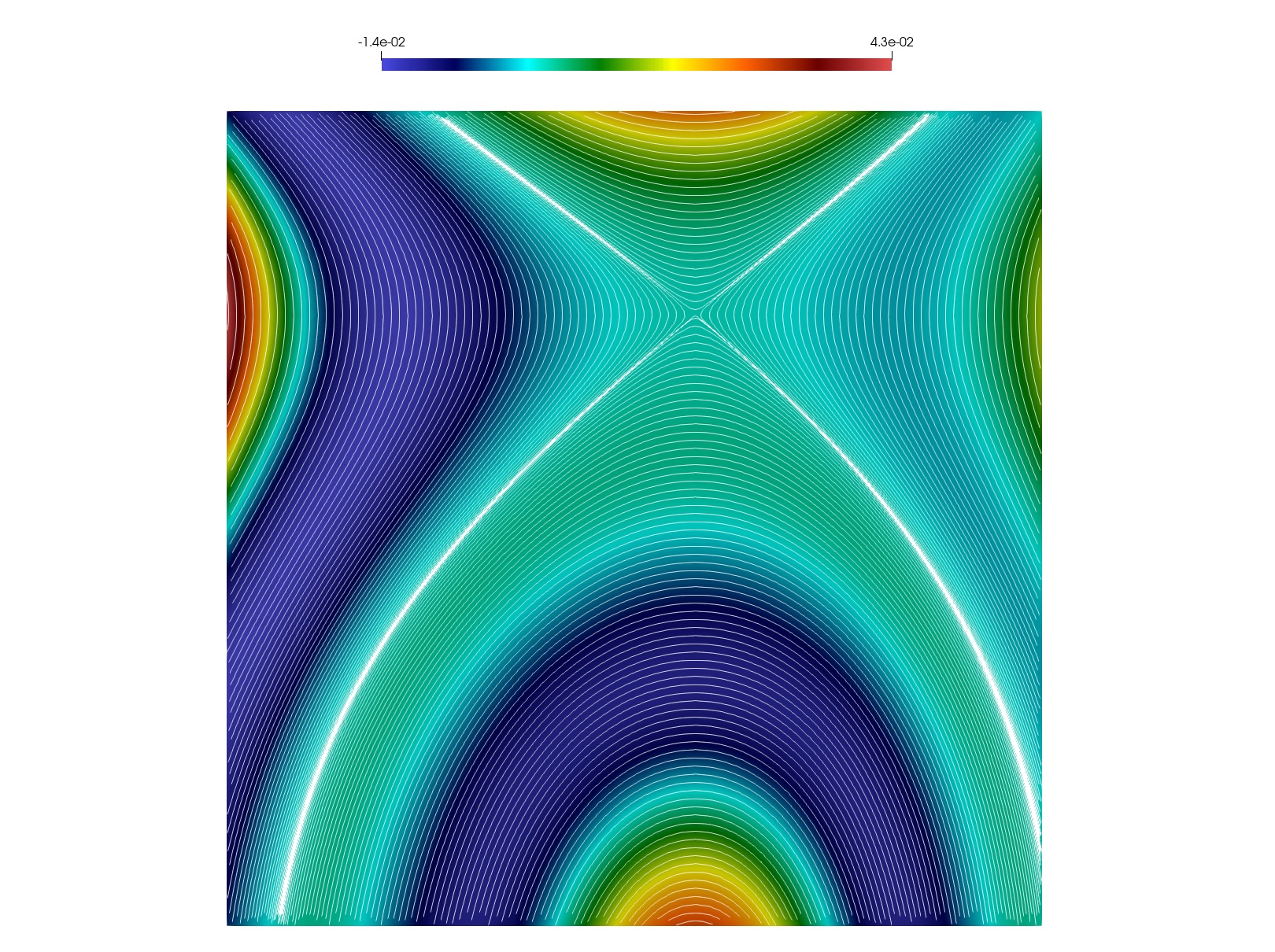} 
\caption{$k_{\parallel}/k_{\perp} = 10^9$}
\end{subfigure}
\caption{Eigenvectors, $\phi_l^{\omega_{15}}$ for $l=1,\ldots,6$}
\end{subfigure}
\caption{Local domain with magnetic field and corresponding eigenvectors for three different anisotropy ratios.}
\label{fig:spectral}
\end{figure}

In order to choose the most ``meaningful'' eigenvectors,  we order eigenvalues  $\{\lambda_j\}$ in increasing magnitude $\lambda_1 \leq \lambda_2 \leq \ldots \leq \lambda_j \leq \ldots$ with corresponding eigenvectors $\{\phi^{\omega_i}_j\}$. The spectral multiscale basis functions with the smallest eigenvalues follow heat flow in the direction of the magnetic field, as represented in Figure \ref{fig:spectral}. We plot the first six eigenvectors in the local domain $\omega_{15}$ on a $5 \times 5$ coarse grid for three anisotropy ratios $k_{\parallel}/k_{\perp} = 10^1, 10^3$ and $10^9$. The first eigenvector is constant or near constant (to remove unnecessary oscillations, we artificially make the first basis constant in the presence of small oscillations). From the next five eigenvectors, we observe how the anisotropy ratio affects them. For small anisotropy, the basis functions largely respect the magnetic field lines in the parallel direction, but the perpendicular diffusion must also be accounted for, and we can see that the basis functions also represent heat flow in that direction. In contrast, for high anisotropy we observe that eigenvectors are more or less exclusively aligned with the magnetic field lines, only representing heat flow in the parallel direction.

\subsection{Coarse grid approximation}

We choose eigenvectors corresponding to the $J^{\omega_i}$ smallest eigenvalues and define a multiscale space
\[
V_H = \text{span}\{ \psi^{\omega_i}_j  = \chi_i \phi^{\omega_i}_j, \ j=1,\ldots,J^{\omega_i}, \ i=1,\ldots,N_H^{vert}\},
\]
where $\chi_i$ is the linear partition of unity functions and $N_H^{vert}$ is the number of local domains (number of coarse grid vertices). 

The multiscale solution $T_{ms} \in V_H$ can be obtained as follows
\begin{equation}
\label{eq:7}
\frac{1}{\tau} m (T^n_{ms} - T^{n-1}_{ms}, v)  + a(T^n_{ms}, v) = l(v), \quad \forall v \in V_H,
\end{equation}
where $T^n_{ms} = \sum_{i=1}^{N_H^{vert}} \sum_{j=1}^{J^{\omega_i}} T^n_{i,j} \ \psi^{\omega_i}_j$. 

From an implementation perspective, we create an interpolation operator from the coarse to the fine grid ($P$) and Galerkin restriction operator from the fine to the coarse grid ($R = P^T$), where 
\begin{equation} 
\label{ms-r}
P = \left[ 
\psi^{\omega_1}_1, \ldots, \psi^{\omega_1}_{J^{\omega_1}}
\ldots
\psi^{\omega_{N_H^{vert}}}_1, \ldots, \psi^{\omega_{N_H^{vert}}}_{J^{\omega_{N_H^{vert}}}}
\right].
\end{equation} 
where the spectral multiscale basis functions $\psi^{\omega_i}_j$ are mapped from local indexing in $\mathcal{T}^{\omega_i}_h$ to global indexing in $\mathcal{T}_h$.  

On the coarse grid we then have the following reduced order implicit equation: 
\begin{equation}
\label{eq:ms}
\frac{1}{\tau} M_H  T^n_H +  A_H T^n_H =  \frac{1}{\tau} M_H  T^{n-1}_H + F_H,
\end{equation}
with 
\[
M_H = P^T M_h P, \quad 
A_H = P^T A_h P, \quad 
F_H = P^T F_h,
\]
where $T^n_H$ is the solution on the coarse grid at time $t_n$. The fine-scale solution is reconstructed by
\[
T_{ms}^n = P \ T^n_H.
\]
The size of the coarse equation \eqref{eq:ms} is $DOF_H = \sum_{i=1}^{N_H^{vert}} J^{\omega_i}$, with $J^{\omega_i}$ denoting the number of local multiscale basis functions in $\omega_i$. In our numerical approach, we take the same number of basis functions in each local domain $\omega_i$ ($J^{\omega_i} = J$), then $DOF_H = J \cdot N_H^{vert}$.

\subsection{Convergence of the multiscale space}

In this section, we analyze the convergence of the multiscale space induced by a generalized spectral problem \eqref{eq:sp} to understand how the number of multiscale basis functions and coarse grid resolution affect the coarse grid approximation error. 
We consider the following norms \cite{abreu2019convergence,efendiev2011multiscale, vasilyeva2024generalized}
\begin{align*}
& ||u||^2_{D_h} = u^T D_h u  = (u, u)_{D_h}, \quad 
||u||^2_{A_h} = u^T A_h u =  (u, u)_{A_h}, \\ 
& ||B_h u||^2_{D_h} 
= (B_h u)^T D_h (B_h u) = (A_h u)^T D^{-1}_h (A_h u)
= (B_h u, B_h u)_{D_h} = (A_h u, A_h u)_{D^{-1}_h} ,
\end{align*}
with  $B_h = D^{-1}_h A_h$, $D_h = D_h^T > 0$ and $A_h = A_h^T > 0$.

To construct spectral multiscale basis functions, we solve a generalized eigenvalue problem \eqref{eq:sp} in each subdomain $\omega_i$ and choose eigenvectors $\phi_j^{\omega_i}$ ($j = 1,\ldots,J^{\omega_i}$) corresponding to the $J^{\omega_i}$ smallest eigenvalues $\lambda_j^{\omega_i}$.  
The eigenvectors $\{\phi^{\omega_i}_j\}$ form an orthonormal basis with respect to the inner product $(u, v)_{D_h^{\omega_i}}$. 
Then, in each local domain $\omega_i$, we define the local projection as follows
\begin{equation}
\label{pl}
P_{J^{\omega_i}}^{\omega_i} v 
= \sum_{j=1}^{J^{\omega_i}} (v, {\phi}_j^{\omega_i})_{D_h^{\omega_i}} {\phi}_j^{\omega_i},
\quad v \in V_h.
\end{equation}
For the local projection the following inequalities hold (see \citep{abreu2019convergence} for details):
\begin{equation}
\label{plest}
\begin{split}
&||v - P_{J^{\omega_i}}^{\omega_i} v ||^2_{D^{\omega_i}_h } 
\leq \frac{1}{ \lambda_{J^{\omega_i}+1}^{\omega_i} } ||v||^2_{A^{\omega_i}_h},
\quad 
||v - P_{J^{\omega_i}}^{\omega_i} v ||^2_{A^{\omega_i}_h} 
\leq \frac{1}{ \lambda_{J^{\omega_i}+1}^{\omega_i} } 
||B_h^{\omega_i} v||^2_{D^{\omega_i}_h},
\\ 
&||v - P_{J^{\omega_i}}^{\omega_i} v ||^2_{D^{\omega_i}_h} 
\leq \frac{1}{ (\lambda_{J^{\omega_i}+1}^{\omega_i})^2 } 
||B_h^{\omega_i} v||^2_{D^{\omega_i}_h}.
\end{split}
\end{equation}
The first inequality corresponds to a local weak approximation property over $\omega_i$ in multigrid literature \cite{vassilevski2008multilevel}. Assuming unit diagonal so that $D^{\omega_i}_h = I$, the second inequality corresponds to a local strong approximation property \cite{vassilevski2008multilevel}, and the third a local fractional approximation property FAP$(1,0)$ \cite[Def. 1]{manteuffel2019convergence} (these latter approximation properties are not generalized in multigrid literature to $D^{\omega_i}_h \neq I$).

Next, we write the coarse interpolation $\Pi$: $V_h \rightarrow V_H$ as follows
\begin{equation}
\label{pg}
\Pi v = \sum_{i=1}^{N_H^{vert}}  \chi_i (P_{J^{\omega_i}}^{\omega_i} v).
\end{equation}
and $v - \Pi v  = \sum_{i=1}^N \chi_i (v - P_{J^{\omega_i}}^{\omega_i}  v)$. 
Using properties of partition of unity functions $\chi_i$ ($\chi_i \leq 1$ and $|\nabla  \chi_i| \leq 1/H^2$) \cite{melenk1996partition}, we can obtain the following estimates for the global projection 
\[
\begin{split}
||u - \Pi u ||_{D_h}^2 
& = \sum_K ||u - \Pi u ||_{D^K_h}^2 \\
& \preceq
\sum_K  \sum_{x_l \in K}  ||\chi_l (u - P_{J^{\omega_l}}^{\omega_l}  u)||^2_{D^{\omega_l}_h}
 \preceq
\sum_K  \sum_{x_l \in K} ||u - P_{J^{\omega_l}}^{\omega_l}  u||_{D^{\omega_l}_h}^2,
\\
||u - \Pi u ||_{A_h}^2
 & \preceq
\sum_K  \sum_{x_l \in K}  ||\chi_l (u - P_{J^{\omega_l}}^{\omega_l}  u)||^2_{A^{\omega_l}_h} \\
 & \preceq
\sum_K  \sum_{x_l \in K}  \frac{1}{H^2} ||u - P_{J^{\omega_l}}^{\omega_l}  u||_{D^{\omega_l}_h}^2  
+  
\sum_K  \sum_{x_l \in K}  ||u - P_{J^{\omega_l}}^{\omega_l}  u||_{A^{\omega_l}_h}^2.
\end{split}
\]

We let  
$\lambda_{J+1} = \min_K  \lambda_{K, J+1}$ and 
$\lambda_{K, J+1} = \min_{x_l \in K}  \lambda_{J^{\omega_l}+1}^{\omega_l}$, therefore by combing with estimates for local projection \eqref{plest} we obtain
\begin{equation}
\begin{split}
\label{pgest}
||v - \Pi v ||_{D_h}^2 & \leq 
\frac{1}{\lambda_{J+1}}  ||v||^2_{A_h},
\\
||v - \Pi v ||_{A_h}^2 & \leq 
\left( \frac{1}{H^2 \lambda_{J+1}^2} + \frac{1}{\lambda_{J+1}} \right) ||B_h v||^2_{D_h},
\\
||v - \Pi v ||_{D_h}^2 & \leq 
\frac{1}{\lambda_{J+1}^2}  ||B_h v||^2_{D_h}. 
\end{split}
\end{equation}
Here we have that $\Pi$ satisfies weak, strong, and fractional FAP$(1,0)$ approximation properties, respectively, with constants as shown above. As discussed in \ref{app:two-grid}, the weak approximation property provides necessary and sufficient conditions for convergence of a two-grid method. 

Finally under some additional regularity and appropriate initial conditions, we have the following error estimate of the multiscale method  (see \ref{app:ms})
\begin{equation}
\label{eq:est}
||T^n_h - T^n_{ms}||_{M_h}^2 + \tau \sum_{k=1}^n || T^k_h - T^k_{ms}||_{A_h}^2
\preceq 
|| T_h^0 - T^0_{ms} ||^2_{M_h}  \\
+ \tau  \sum_{k=1}^n \frac{H^2}{\Lambda^*}   || B_h T_h^k||^2_{D_h},
\end{equation}
where $\Lambda^* = \lambda_{J+1} H^2$ (local domain is scaled to a domain of size one),  $T_h^n \in V_h$ and $T^n_{ms} \in V_H$ are the fine-scale and multiscale solutions from \eqref{eq:5} and \eqref{eq:7}, respectively.


\section{Two-grid multiscale preconditioner}\label{sec:tg}

Next, we consider the construction of the iterative solver, where we use the constructed multiscale space as a two-grid preconditioner that converges independently of the contrast of the parallel and perpendicular heat conductivities. Such contrast is challenging for a typical multigrid preconditioner. In contrast with classic geometric or algebraic restriction and interpolation operators, our transfer operators built on local spectral problems lead to a very accurate coarse approximation, as shown above. 

Following equation \eqref{eq:6}, for each implicit time step we solve a system of linear equations
\begin{equation}
\label{eq:sf}   
Q_h T^n_h = b^{n-1}_h,
\end{equation}
where 
\[
Q_h = \frac{1}{\tau} M_h + A_h, \quad 
b^{n-1}_h = \frac{1}{\tau} M_h T^{n-1}_h + F_h.
\]

Following equation \eqref{eq:ms}, we have the following system on the coarse grid
\begin{equation}
\label{eq:sc}  
Q_H T^n_H = b^{n-1}_H, \quad T^n_{ms} = P \ T^n_H,
\end{equation}
where 
\[
Q_H = \frac{1}{\tau} M_H + A_H, \quad 
b^{n-1}_H = \frac{1}{\tau} M_H  T^{n-1}_H + F_H.
\]

Next, we follow the framework of the Algebraic Multigrid Method (AMG) method in constructing a two-grid solver \cite{ruge1987algebraic, xu2017algebraic, vassilevski2008multilevel, falgout2005two, vassilevski2011coarse}.
For a given initial guess $y_h^{(0)} = T^{n-1}_h$ in the two-grid algorithm and smoothing operator $S\approx Q_h$, we have:
\begin{enumerate}
\item \textit{Pre-smoothing}: 
\[
y_h^{(1)} =  y_h^{(0)} + S^{-1} r_h^{(0)}, \quad 
r_h^{(0)} = (b^{n-1}_h - Q_h y_h^{(0)})
\]
\item \textit{Coarse-grid correction}:
\begin{enumerate}
\item[2.1]   \textit{Restriction}:  
\[
r_H = P^T r_h^{(1)}, \quad 
r_h^{(1)} = (b^{n-1}_h - Q_h y_h^{(1)}).
\]
\item[2.2]  \textit{Coarse-grid solution}: 
\[
Q_H e_H = r_H.
\]
\item[2.3] \textit{Interpolation and update}: 
\[
y^{(2)} = y^{(1)} + P e_H.
\]
\end{enumerate}
\item \textit{Post-smoothing}: 
\[
y_{TG} =  y^{(2)} + S^{-T} r_h^{(2)}, \quad 
r_h^{(2)} = (b^{n-1}_h - Q_h y_h^{(2)}).
\]
\end{enumerate}

In the two-grid algorithm, we have the following error transfer operator \cite{falgout2005two,notay2007convergence}
\[
E_{TG}  = (I - S^{-T} Q_h) (I - P Q_H^{-1} P^T Q_h) (I - S^{-1} Q_h).
\]
In this work, we use five pointwise Jacobi and symmetric Gauss-Seidel iteration pre- and post-smoothing iterations  to remove high-frequency errors.
A multiscale coarse-grid correction is used to attenuate the remaining error not effectively reduced by smoothing.  

The two-grid error propagation can be expressed as follows \cite{falgout2005two, vassilevski2008multilevel}
\[
E_{TG} = I - C^{-1}_{TG} Q_h,
\]
with \[
C^{-1}_{TG} = \bar{S}^{-1} + (I - S^{-T} Q_h) P Q^{-1}_H P^T (I - Q_h S^{-1}),
\]
where $\bar{S} = S (S + S^T  - Q_h)^{-1} S^T$ is the so-called symmetric smoother. 

For the method's convergence, we have  \cite{vassilevski2008multilevel, brezina2011smoothed, falgout2005two}
\[
0 \leq v^T Q_h E_{TG} v \leq \left( 1 - \frac{1}{K_{TG}} \right) v^T Q_h v,
\] 
with
\[
v^T Q_h v \leq v^T C_{TG} v \leq K_{TG} v^T Q_h v,
 \quad \textnormal{where }
 K_{TG} \coloneqq \text{Cond} (C_{TG}^{-1} Q_h).
\]

Using properties of the Jacobi and Gauss-Seidel smoothers, spectral equivalence to the diagonal part of the matrix $Q_h$, and approximation properties of the spectral multiscale space discussed in the previous section, we obtain (see \ref{app:two-grid} for details) 
\[
K_{TG} = \left(1 + \frac{C}{\tau} \right) \frac{H^2}{\Lambda^*}.
\]
Therefore, given a sufficient number of multiscale basis functions, we can obtain anisotropy-independent convergence of the iterative solver. More importantly, in Section \ref{sec:results:mg} we show that $\mathcal{O}(1)$ local basis functions are sufficient for robust anisotropy-independent convergence.

\section{Numerical results}\label{sec:results}

We consider anisotropic heat flow in a square domain $\Omega = [0,1]^2$, with $\kappa_\perp = 1$ and $\kappa_\| \gg 1$ such that the ratio $\kappa_\perp/\kappa_\|$ ranges from $10^3$ to $10^{12}$. We simulate with $t_{max} = 5 \cdot 10^{-6}$ using 10 time steps. Note, the time step is fixed for all anisotropies, which makes for a relatively small timestep at the smallest anisotropy, but a large time step and very stiff system for the largest anisotropy ratio. Although the focus of this paper is on spatial discretization, for realistically high anisotropy ratios a steady formulation of the problem is extremely ill conditioned/bordering on ill-posed, hence why we consider a short time evolution formulation to test the multiscale spatial representation and preconditioner.
The fine grid contains 
100,292 triangular cells, 150,838 facets and 50,547 vertices. 
The fine grid solution using a CG FEM approximation with $P^2$ elements is used as a reference solution.  
Then the number of unknowns of the fine grid problem is 
$DOF_h = N_h = N_h^{vert} + N_h^{e} = 201,385$, 
where $N_h^{vert}$ is the number of fine grid vertices and $N_h^{e}$ is the number of fine grid facets.

\begin{figure}[h!]
\centering
\begin{subfigure}[b]{0.32\textwidth}
\centering
\includegraphics[width=1 \textwidth]{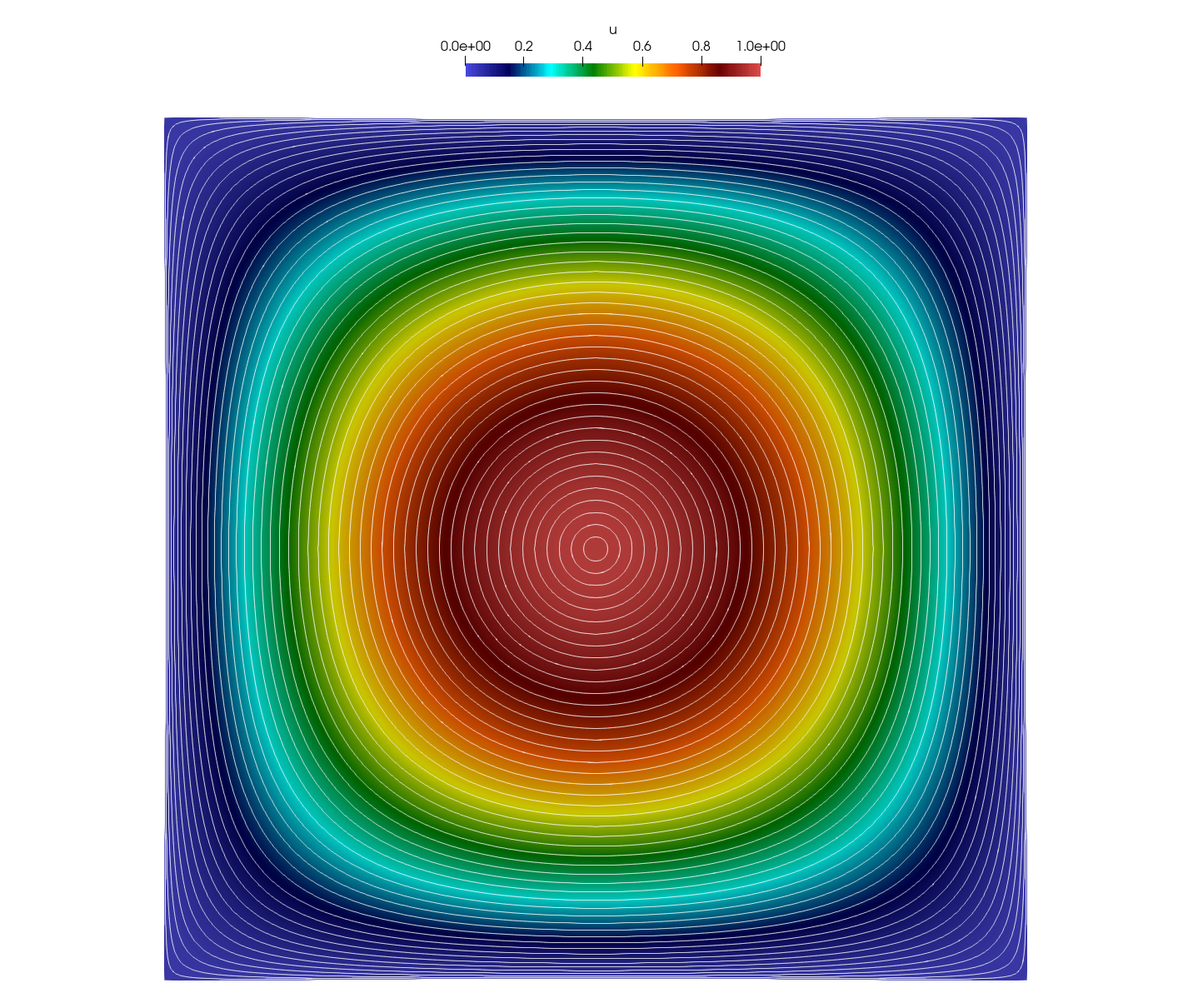} 
\caption{Test 1}
\end{subfigure}
\begin{subfigure}[b]{0.32\textwidth}
\centering
\includegraphics[width=1 \textwidth]{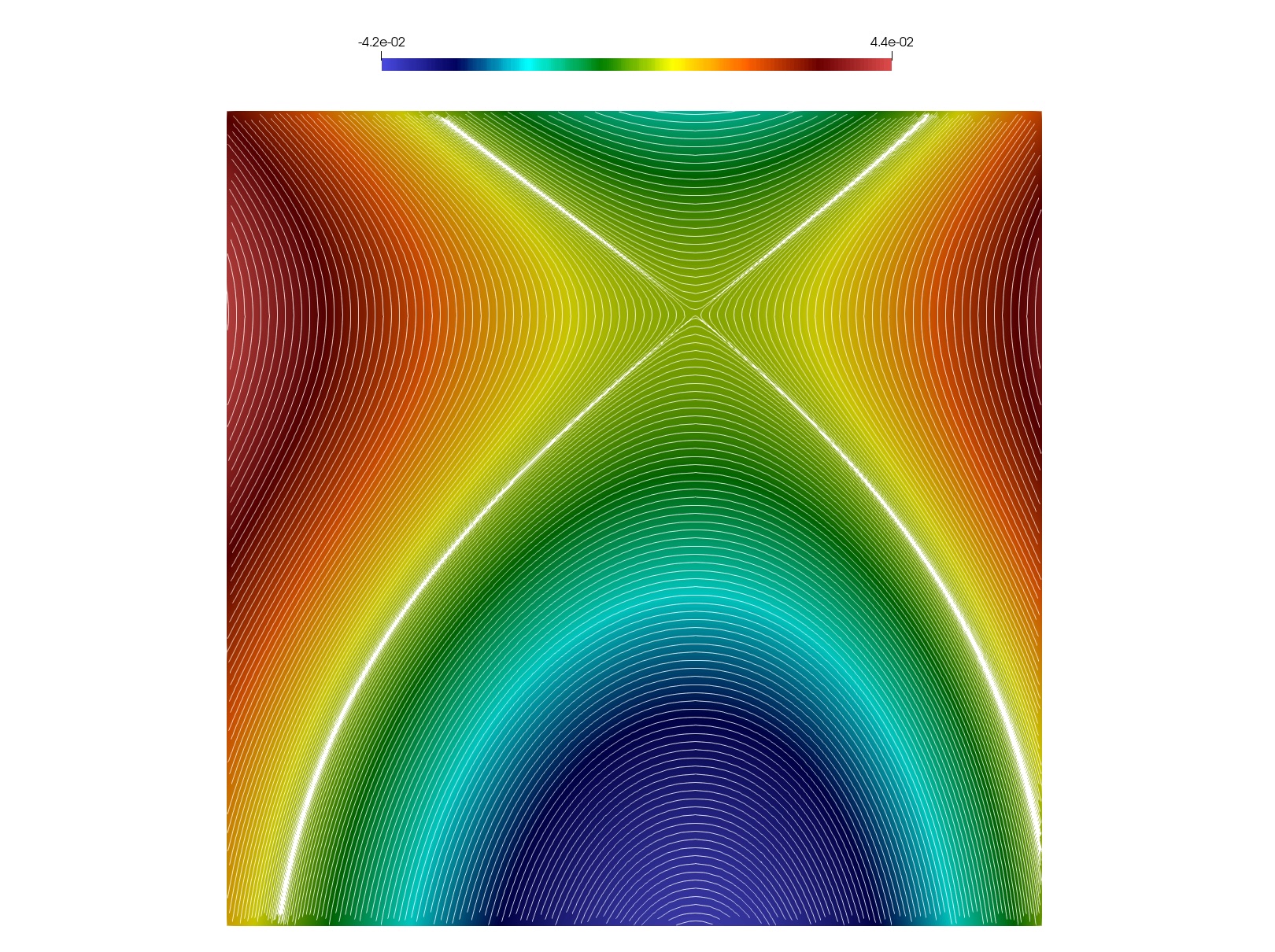} 
\caption{Test 2}
\end{subfigure}
\begin{subfigure}[b]{0.32\textwidth}
\centering
\includegraphics[width=1 \textwidth]{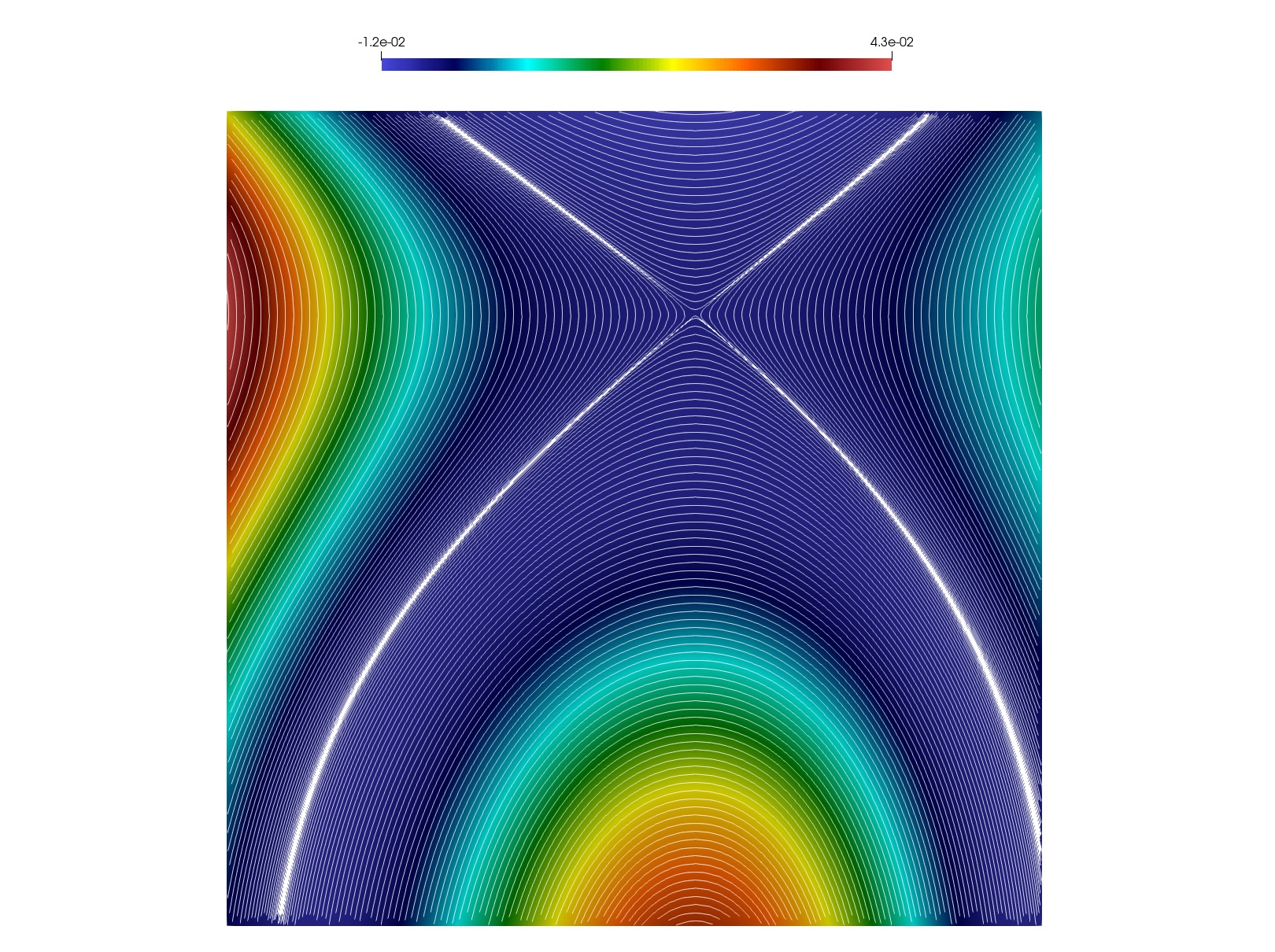} 
\caption{Test 3}
\end{subfigure}
\caption{Reference solutions with magnetic field lines at final time for Test 1, 2 and 3,  $DOF_h=201,385$}
\label{fig:u}
\end{figure}

We simulate three test cases represented in Figure \ref{fig:u}. Each case is characterized by a given magnetic field (Test 1, 2, and 3). In Figure \ref{fig:u}, we plot the temperature distribution with magnetic field lines depicted in white lines. 
We set a source term $f = \kappa_\perp \Delta T_0$ that simulates a counter-forcing. 
The implementation is based on the FEniCS library \cite{logg2012automated}. To solve the spectral problem for the multiscale basis functions construction, we use the sparse package from scipy \cite{virtanen2020scipy}. The Gmsh library \cite{geuzaine2009gmsh} is used to construct a fine-grid conformed with coarse grid edges.

\subsection{Multiscale space approximation}\label{sec:results:ms}

To investigate the accuracy of the proposed multiscale method, we consider  $10 \times 10$ and $20 \times 20$ coarse grids. 
The size of the resulting coarse-grid system is $DOF_H = \sum_{i=1}^{N_H^{vert}} J^{\omega_i}$, with $J^{\omega_i}$ denoting the number of local multiscale basis functions in $\omega_i$. In our numerical approach, we take the same number of basis functions in each local domain $\omega_i$ ($J^{\omega_i} = J$), then $DOF_H = J \cdot N_H^{vert}$. Then we have $DOF_H = J \cdot 121$ for the $10 \times 10$ coarse grid and $DOF_H = J \cdot 441$ for the $20 \times 20$ coarse grid. We use a direct solver with $LU$ factorization for coarse-scale approximations.

\begin{table}[h!]
\centering
\begin{tabular}{|c|c|cc|cc|cc|cc|}
\hline
\multirow{ 2}{*}{$J$}  & \multirow{ 2}{*}{$DOF_H$} 
& $err$ & tm$_{s}$
& $err$ & tm$_{s}$
& $err$ & tm$_{s}$
& $err$ & tm$_{s}$\\ 
  &
  & \multicolumn{2}{c|}{$10^3$} 
  & \multicolumn{2}{c|}{$10^6$}  
  & \multicolumn{2}{c|}{$10^9$}   
  & \multicolumn{2}{c|}{$10^{12}$}\\ \hline
 \multicolumn{10}{|c|}{Coarse grid, $10 \times 10$} \\ \hline  
1  & 121 & 5.53e-03 & 0.020 & 2.77e-01 & 0.020 & 1.00e+00 & 0.020 & 1.00e+00 & 0.020 \\ \hline
2  & 242 & 6.59e-03 & 0.036 & 8.25e-02 & 0.037 & 2.34e-01 & 0.037 & 1.00e+00 & 0.037 \\ \hline
3  & 363 & 2.84e-03 & 0.053 & 3.02e-02 & 0.054 & 5.22e-02 & 0.054 & 9.02e-01 & 0.054 \\ \hline
4  & 484 & 1.53e-03 & 0.071 & 4.58e-03 & 0.073 & 6.28e-03 & 0.073 & 3.43e-01 & 0.072 \\ \hline
6  & 726 & 6.57e-04 & 0.106 & 2.76e-03 & 0.106 & 9.70e-04 & 0.105 & 1.00e-02 & 0.106 \\ \hline
8  & 968 & 2.78e-04 & 0.140 & 3.11e-04 & 0.140 & 2.69e-04 & 0.141 & 4.47e-03 & 0.143 \\ \hline
12  & 1452 & 1.20e-04 & 0.211 & 3.07e-05 & 0.210 & 2.08e-05 & 0.211 & 4.26e-04 & 0.213 \\ \hline
16  & 1936 & 5.72e-05 & 0.283 & 4.15e-06 & 0.281 & 9.49e-06 & 0.281 & 2.08e-04 & 0.284 \\ \hline
24  & 2904 & 3.56e-05 & 0.426 & 4.51e-07 & 0.427 & 2.81e-06 & 0.427 & 4.79e-04 & 0.434 \\ \hline
32  & 3872 & 2.97e-05 & 0.576 & 3.96e-07 & 0.573 & 2.68e-06 & 0.576 & 3.69e-04 & 0.579 \\ \hline
40  & 4840 & 1.87e-05 & 0.725 & 4.18e-07 & 0.726 & 2.61e-06 & 0.730 & 3.82e-04 & 0.731 \\ \hline
48  & 5808 & 1.37e-05 & 0.890 & 4.15e-07 & 0.887 & 2.40e-06 & 0.892 & 2.32e-04 & 0.892 \\ \hline
56  & 6776 & 9.48e-06 & 1.048 & 3.95e-07 & 1.049 & 2.20e-06 & 1.055 & 2.11e-04 & 1.049 \\ \hline
64  & 7744 & 6.24e-06 & 1.254 & 3.25e-07 & 1.219 & 2.16e-06 & 1.217 & 1.90e-04 & 1.218 \\ \hline
 \multicolumn{10}{|c|}{Coarse grid, $20 \times 20$} \\ 
\hline  
1  & 441 & 1.48e-03 & 0.018 & 1.23e-01 & 0.018 & 1.00e+00 & 0.019 & 1.00e+00 & 0.018 \\ \hline
2  & 882 & 2.12e-03 & 0.033 & 4.21e-02 & 0.034 & 4.64e-01 & 0.033 & 1.00e+00 & 0.034 \\ \hline
3  & 1323 & 1.24e-03 & 0.050 & 2.18e-02 & 0.049 & 4.10e-02 & 0.049 & 8.91e-01 & 0.050 \\ \hline
4  & 1764 & 5.09e-04 & 0.065 & 8.45e-04 & 0.065 & 1.78e-03 & 0.065 & 5.32e-02 & 0.066 \\ \hline
6  & 2646 & 3.57e-04 & 0.098 & 1.72e-04 & 0.098 & 1.07e-04 & 0.098 & 6.17e-03 & 0.099 \\ \hline
8  & 3528 & 1.14e-04 & 0.131 & 2.64e-05 & 0.132 & 2.07e-05 & 0.131 & 1.46e-03 & 0.133 \\ \hline
12  & 5292 & 3.42e-05 & 0.203 & 1.09e-06 & 0.200 & 2.09e-06 & 0.203 & 2.80e-04 & 0.203 \\ \hline
16  & 7056 & 1.72e-05 & 0.278 & 3.84e-07 & 0.278 & 1.94e-06 & 0.280 & 2.15e-04 & 0.281 \\ \hline
24  & 10584 & 5.51e-06 & 0.446 & 2.68e-07 & 0.446 & 1.48e-06 & 0.446 & 1.72e-04 & 0.447 \\ \hline
32  & 14112 & 3.63e-06 & 0.641 & 3.31e-08 & 0.639 & 5.29e-07 & 0.642 & 5.85e-05 & 0.644 \\ \hline
40  & 17640 & 1.88e-06 & 0.856 & 3.16e-08 & 0.854 & 5.35e-07 & 0.860 & 2.62e-04 & 0.858 \\ \hline
48  & 21168 & 1.16e-06 & 1.108 & 2.20e-08 & 1.089 & 3.96e-07 & 1.101 & 5.52e-05 & 1.096 \\ \hline
56  & 24696 & 7.97e-07 & 1.361 & 1.90e-08 & 1.361 & 5.09e-07 & 1.368 & 1.03e-04 & 1.352 \\ \hline
64  & 28224 & 5.91e-07 & 1.648 & 4.32e-09 & 1.636 & 5.96e-07 & 1.665 & 1.28e-04 & 1.657 \\ \hline
\end{tabular}
\caption{Test 1. Relative L$_2$ error ($err$) with time of solution (tm$_{s}$)}
\label{table-err1}
\end{table}

\begin{table}[h!]
\centering
\begin{tabular}{|c|c|cc|cc|cc|cc|}
\hline
\multirow{ 2}{*}{$J$}  & \multirow{ 2}{*}{$DOF_H$} 
& $err$ & tm$_{s}$
& $err$ & tm$_{s}$
& $err$ & tm$_{s}$
& $err$ & tm$_{s}$\\ 
  &
  & \multicolumn{2}{c|}{$10^3$} 
  & \multicolumn{2}{c|}{$10^6$}  
  & \multicolumn{2}{c|}{$10^9$}   
  & \multicolumn{2}{c|}{$10^{12}$}\\ \hline
 \multicolumn{10}{|c|}{Coarse grid, $10 \times 10$} \\ \hline  
1  & 121 & 1.24e-02 & 0.019 & 5.57e-01 & 0.020 & 1.00e+00 & 0.019 & 1.00e+00 & 0.020 \\ \hline
2  & 242 & 1.02e-02 & 0.037 & 1.08e-01 & 0.037 & 9.91e-01 & 0.038 & 1.00e+00 & 0.037 \\ \hline
3  & 363 & 5.03e-03 & 0.054 & 4.23e-02 & 0.056 & 5.78e-02 & 0.054 & 1.00e+00 & 0.055 \\ \hline
4  & 484 & 2.89e-03 & 0.072 & 5.58e-03 & 0.072 & 2.85e-02 & 0.072 & 9.93e-01 & 0.072 \\ \hline
6  & 726 & 1.55e-03 & 0.107 & 1.54e-03 & 0.106 & 7.71e-04 & 0.107 & 1.51e-01 & 0.106 \\ \hline
8  & 968 & 8.27e-04 & 0.143 & 4.23e-04 & 0.141 & 4.39e-04 & 0.142 & 4.94e-02 & 0.141 \\ \hline
12  & 1452 & 3.84e-04 & 0.212 & 2.82e-05 & 0.210 & 1.11e-04 & 0.213 & 1.32e-02 & 0.212 \\ \hline
16  & 1936 & 1.92e-04 & 0.283 & 9.63e-06 & 0.282 & 6.07e-05 & 0.283 & 7.17e-03 & 0.283 \\ \hline
24  & 2904 & 1.07e-04 & 0.431 & 2.80e-06 & 0.427 & 2.80e-05 & 0.430 & 2.74e-03 & 0.424 \\ \hline
32  & 3872 & 6.73e-05 & 0.579 & 1.90e-06 & 0.574 & 2.46e-05 & 0.579 & 2.61e-03 & 0.577 \\ \hline
40  & 4840 & 3.58e-05 & 0.727 & 1.97e-06 & 0.728 & 2.28e-05 & 0.727 & 2.41e-03 & 0.729 \\ \hline
48  & 5808 & 2.31e-05 & 0.886 & 1.93e-06 & 0.885 & 2.17e-05 & 0.888 & 2.27e-03 & 0.894 \\ \hline
56  & 6776 & 1.40e-05 & 1.053 & 1.87e-06 & 1.055 & 2.01e-05 & 1.055 & 2.09e-03 & 1.059 \\ \hline
64  & 7744 & 9.45e-06 & 1.224 & 1.57e-06 & 1.226 & 1.57e-05 & 1.225 & 1.60e-03 & 1.228 \\ \hline
 \multicolumn{10}{|c|}{Coarse grid, $20 \times 20$} \\ \hline  
1  & 441 & 3.14e-03 & 0.018 & 2.29e-01 & 0.019 & 1.00e+00 & 0.019 & 1.00e+00 & 0.018 \\ \hline
2  & 882 & 3.64e-03 & 0.033 & 6.95e-02 & 0.034 & 8.84e-01 & 0.035 & 1.00e+00 & 0.034 \\ \hline
3  & 1323 & 2.11e-03 & 0.049 & 3.37e-02 & 0.049 & 2.37e-02 & 0.050 & 9.96e-01 & 0.050 \\ \hline
4  & 1764 & 8.88e-04 & 0.066 & 1.64e-03 & 0.067 & 2.65e-03 & 0.069 & 4.75e-01 & 0.065 \\ \hline
6  & 2646 & 5.11e-04 & 0.097 & 1.19e-04 & 0.098 & 3.31e-04 & 0.100 & 2.74e-02 & 0.099 \\ \hline
8  & 3528 & 1.84e-04 & 0.132 & 2.45e-05 & 0.133 & 4.98e-05 & 0.132 & 8.20e-03 & 0.131 \\ \hline
12  & 5292 & 5.65e-05 & 0.204 & 3.66e-06 & 0.202 & 2.77e-05 & 0.203 & 3.60e-03 & 0.202 \\ \hline
16  & 7056 & 3.15e-05 & 0.280 & 1.73e-06 & 0.279 & 1.33e-05 & 0.278 & 1.23e-03 & 0.278 \\ \hline
24  & 10584 & 8.08e-06 & 0.450 & 9.29e-07 & 0.448 & 9.36e-06 & 0.447 & 6.50e-04 & 0.445 \\ \hline
32  & 14112 & 5.38e-06 & 0.642 & 1.82e-06 & 0.640 & 4.70e-06 & 0.641 & 8.35e-04 & 0.638 \\ \hline
40  & 17640 & 3.03e-06 & 0.862 & 1.79e-06 & 0.861 & 5.09e-06 & 0.855 & 8.41e-04 & 0.866 \\ \hline
48  & 21168 & 1.87e-06 & 1.099 & 1.39e-06 & 1.098 & 3.70e-06 & 1.097 & 5.24e-04 & 1.110 \\ \hline
56  & 24696 & 1.29e-06 & 1.378 & 1.26e-06 & 1.365 & 3.81e-06 & 1.372 & 3.03e-04 & 1.378 \\ \hline
64  & 28224 & 8.91e-07 & 1.671 & 2.49e-07 & 1.673 & 1.05e-06 & 1.662 & 1.24e-04 & 1.646 \\ \hline
\end{tabular}
\caption{Test 2. Relative L$_2$ error ($err$) with time of solution (tm$_{s}$)}
\label{table-err2}
\end{table}

\begin{table}[h!]
\centering
\begin{tabular}{|c|c|cc|cc|cc|cc|}
\hline
\multirow{ 2}{*}{$J$}  & \multirow{ 2}{*}{$DOF_H$} 
& $err$ & tm$_{s}$
& $err$ & tm$_{s}$
& $err$ & tm$_{s}$
& $err$ & tm$_{s}$\\ 
  &
  & \multicolumn{2}{c|}{$10^3$} 
  & \multicolumn{2}{c|}{$10^6$}  
  & \multicolumn{2}{c|}{$10^9$}   
  & \multicolumn{2}{c|}{$10^{12}$}\\ \hline
 \multicolumn{10}{|c|}{Coarse grid, $10 \times 10$} \\ \hline  
1  & 121 & 2.26e-03 & 0.019 & 8.82e-03 & 0.020 & 8.78e-03 & 0.020 & 8.51e-03 & 0.020 \\ \hline
2  & 242 & 5.22e-04 & 0.037 & 1.37e-03 & 0.037 & 9.00e-04 & 0.037 & 8.93e-04 & 0.037 \\ \hline
3  & 363 & 4.08e-04 & 0.055 & 2.95e-03 & 0.055 & 6.60e-04 & 0.054 & 3.89e-04 & 0.056 \\ \hline
4  & 484 & 2.35e-04 & 0.071 & 3.23e-04 & 0.072 & 2.17e-04 & 0.072 & 2.16e-04 & 0.071 \\ \hline
6  & 726 & 1.34e-04 & 0.106 & 6.10e-05 & 0.106 & 5.99e-05 & 0.106 & 3.99e-04 & 0.106 \\ \hline
8  & 968 & 8.54e-05 & 0.141 & 1.76e-05 & 0.141 & 6.01e-05 & 0.140 & 1.74e-04 & 0.141 \\ \hline
12  & 1452 & 4.53e-05 & 0.213 & 4.55e-06 & 0.212 & 9.58e-06 & 0.212 & 1.44e-04 & 0.212 \\ \hline
16  & 1936 & 2.33e-05 & 0.282 & 1.56e-06 & 0.285 & 2.14e-06 & 0.284 & 3.66e-05 & 0.284 \\ \hline
24  & 2904 & 9.60e-06 & 0.431 & 6.75e-07 & 0.430 & 8.35e-07 & 0.429 & 6.00e-06 & 0.430 \\ \hline
32  & 3872 & 3.77e-06 & 0.579 & 4.12e-07 & 0.581 & 4.19e-07 & 0.579 & 1.35e-06 & 0.582 \\ \hline
40  & 4840 & 1.96e-06 & 0.734 & 3.32e-07 & 0.732 & 3.73e-07 & 0.733 & 2.34e-06 & 0.726 \\ \hline
48  & 5808 & 1.44e-06 & 0.894 & 3.10e-07 & 0.893 & 3.11e-07 & 0.886 & 4.14e-07 & 0.892 \\ \hline
56  & 6776 & 9.65e-07 & 1.058 & 2.97e-07 & 1.050 & 2.97e-07 & 1.054 & 3.95e-07 & 1.050 \\ \hline
64  & 7744 & 7.56e-07 & 1.219 & 2.89e-07 & 1.235 & 2.89e-07 & 1.225 & 3.02e-07 & 1.231 \\ \hline
 \multicolumn{10}{|c|}{Coarse grid, $20 \times 20$} \\ \hline  
1  & 441 & 5.74e-04 & 0.018 & 5.23e-03 & 0.018 & 5.18e-03 & 0.018 & 4.88e-03 & 0.018 \\ \hline
2  & 882 & 2.04e-04 & 0.033 & 5.58e-04 & 0.034 & 4.86e-04 & 0.034 & 2.24e-04 & 0.034 \\ \hline
3  & 1323 & 1.82e-04 & 0.050 & 1.48e-03 & 0.049 & 3.14e-04 & 0.049 & 7.41e-05 & 0.049 \\ \hline
4  & 1764 & 5.48e-05 & 0.067 & 3.86e-05 & 0.066 & 2.03e-04 & 0.066 & 1.55e-04 & 0.066 \\ \hline
6  & 2646 & 3.36e-05 & 0.099 & 4.81e-06 & 0.098 & 3.87e-05 & 0.098 & 2.11e-04 & 0.100 \\ \hline
8  & 3528 & 1.64e-05 & 0.131 & 1.43e-06 & 0.131 & 3.46e-06 & 0.132 & 1.14e-04 & 0.132 \\ \hline
12  & 5292 & 4.46e-06 & 0.203 & 4.90e-07 & 0.204 & 1.14e-06 & 0.206 & 2.06e-05 & 0.206 \\ \hline
16  & 7056 & 1.89e-06 & 0.281 & 3.23e-07 & 0.282 & 4.13e-07 & 0.280 & 4.33e-06 & 0.281 \\ \hline
24  & 10584 & 6.26e-07 & 0.449 & 2.42e-07 & 0.452 & 2.95e-07 & 0.446 & 1.56e-06 & 0.447 \\ \hline
32  & 14112 & 3.63e-07 & 0.644 & 2.31e-07 & 0.655 & 2.83e-07 & 0.651 & 2.17e-06 & 0.644 \\ \hline
40  & 17640 & 2.25e-07 & 0.870 & 2.13e-07 & 0.861 & 2.31e-07 & 0.864 & 5.63e-07 & 0.858 \\ \hline
48  & 21168 & 1.77e-07 & 1.111 & 2.06e-07 & 1.108 & 2.44e-07 & 1.100 & 2.95e-07 & 1.094 \\ \hline
56  & 24696 & 1.56e-07 & 1.369 & 1.96e-07 & 1.379 & 2.57e-07 & 1.398 & 4.37e-07 & 1.368 \\ \hline
64  & 28224 & 1.44e-07 & 1.658 & 1.82e-07 & 1.713 & 5.02e-07 & 1.703 & 1.87e-06 & 1.646 \\ \hline
\end{tabular}
\caption{Test 3. Relative L$_2$ error ($err$) with time of solution (tm$_{s}$)}
\label{table-err4}
\end{table}

\begin{figure}[h!]
\centering
\includegraphics[width=0.32 \textwidth]{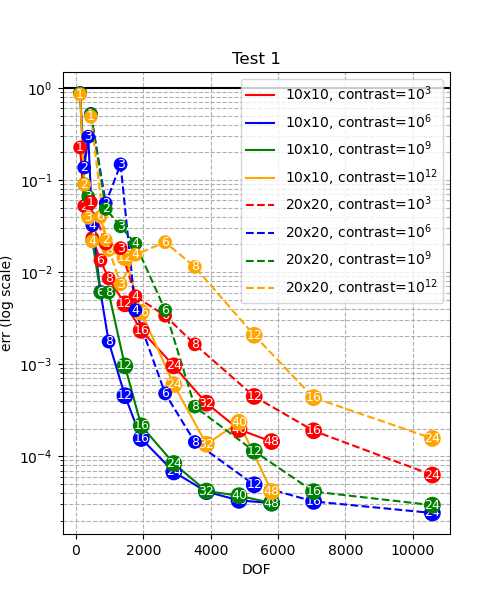} 
\includegraphics[width=0.32 \textwidth]{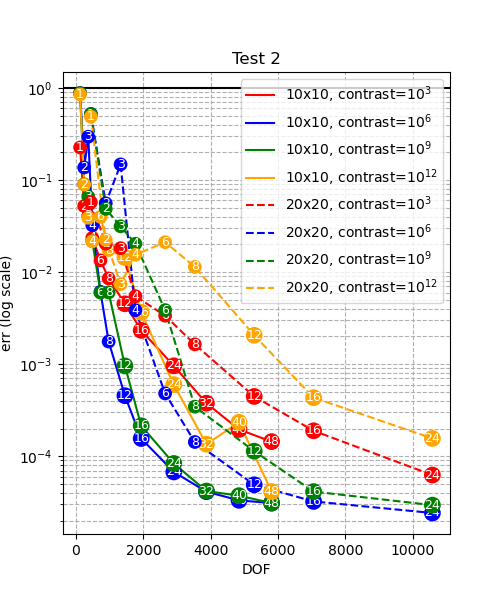} 
\includegraphics[width=0.32 \textwidth]{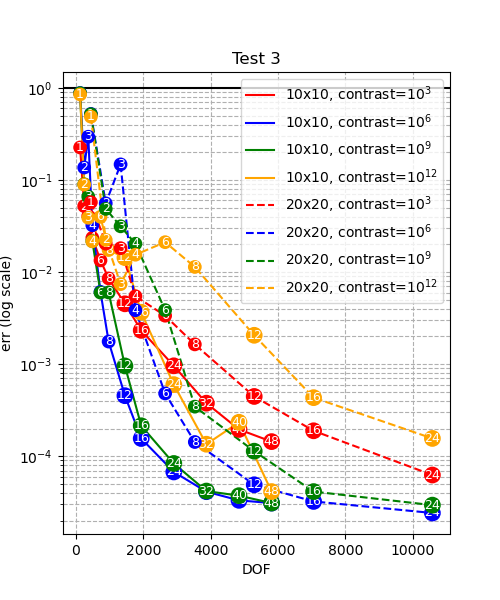}
\caption{Error vs $DOF_H$ for Tests 1, 2 and 3 (from left to right). }
\label{fig:err}
\end{figure}

For error calculations, we use the relative L$_2$ error for the temperature field at the final time
\[
err = \frac{||T_h - T_{ms}||_{L_2}}{||T_h||_{L_2}},
\]
where 
$||u||^2_{L_2} = u^T u = (u, u)$ .

Tables \ref{table-err1},  \ref{table-err2}, and \ref{table-err4} contain numerical results for the presented multiscale method for  Tests 1, 2, and 3, respectively. We investigate the influence of increasing the anisotropy ratio $k_{\parallel}/k_{\perp}$ on the method's robustness, and also present errors within an increasing number of multiscale basis functions $J$, corresponding to the size of the coarse grid system $DOF_H$. The numerical results are presented for $10 \times 10$ and $20 \times 20$ coarse grids. The solution time is measured without taking into account the time it takes to construct the multiscale basis functions (restriction operator), which can be done in the offline stage. Therefore, the time indicated (denoted as tm$_s$) refers to the time required for online computations. It is worth noting that the offline stage of constructing the multiscale basis can be carried out in a fully parallel manner. Additionally, the computation time of the basis in a local domain depends on the size of the coarse grid, so it is important to find a balance between the coarse grid size and fine grid resolution in order to achieve faster calculations in the local domain.

Figure \ref{fig:err} illustrates the results presented in Tables \ref{table-err1},  \ref{table-err2}, and \ref{table-err4}. The graph shows how the relative L$_2$ error decreases with respect to the size of the coarse grid system, $DOF_H$, to demonstrate the relationship between the number of basis functions $M$ and the coarse grid resolution $H$.
Generally, the error increases as $J$ and coarse grid resolution increases in all test cases with four contrasts. The number of basis functions directly affects the system size, and therefore, the solution time is increased proportionally.  
We observe that higher anisotropy requires a larger number of basis functions to capture the fine-scale behavior. 
We also see that on the $20 \times 20$ coarse grid, we can obtain good results with a smaller number of multiscale basis functions. However, when comparing accuracy for a fixed number of coarse system DOFs, better results are typically achieved using the $10\times 10$ grid with more multiscale basis functions. 
%
In Test 2, the coarse grid $10 \times 10$ and $20 \times 20$ results show a similar trend as in Test 1. 
We notice that for Tests 1 and 2 the proposed method gives a larger errors compared to Test 3. In Test 3, we achieve good results with a small number of basis functions $J$. Additionally, the impact of the coarse grid size is not as significant as in Tests 1 and 2. 
In all cases, we observe that we can achieve high levels of accuracy with a careful choice of coarse grid resolution and number of multiscale basis functions to balance the efficiency of computation and construction.

\subsection{Two-grid multiscale preconditioner}\label{sec:results:mg}

%
In this section, we consider the performance of the proposed multiscale solver as a two-grid preconditioner for the conjugate gradient iterative solver from Scipy sparse library \cite{virtanen2020scipy}. As a smoother, we use five Jacobi and Gauss-Seidel iterations ($\nu=5$) from pyamg library \cite{bell2022pyamg, bell2023pyamg}. To stop iterations, we use a default set of parameters in PCG with relative tolerance $r_{tol}= 10^{-5}$ and set 100 as the maximum number of iterations (nc indicates no convergence within the given restrictions). 

\begin{table}[h!]
\centering
\begin{tabular}{|c||cc|cc|cc|cc||cc|cc|cc|cc|}
\hline
\multirow{ 2}{*}{$J$} 
 & $\bar{N}$ & tm$_{s}$
 & $\bar{N}$ & tm$_{s}$
 & $\bar{N}$ & tm$_{s}$
 & $\bar{N}$ & tm$_{s}$
 & $\bar{N}$ & tm$_{s}$
 & $\bar{N}$ & tm$_{s}$
 & $\bar{N}$ & tm$_{s}$
 & $\bar{N}$ & tm$_{s}$\\ 
  & \multicolumn{2}{c|}{$10^3$} 
  & \multicolumn{2}{c|}{$10^6$}  
  & \multicolumn{2}{c|}{$10^9$}   
  & \multicolumn{2}{c||}{$10^{12}$}
  & \multicolumn{2}{c|}{$10^3$} 
  & \multicolumn{2}{c|}{$10^6$}  
  & \multicolumn{2}{c|}{$10^9$}   
  & \multicolumn{2}{c|}{$10^{12}$}\\ \hline
&\multicolumn{8}{|c||}{Jacobi smoother}&\multicolumn{8}{|c|}{Gauss-Seidel smoother} \\ 
\hline
& \multicolumn{16}{|c|}{Coarse grid, $10 \times 10$} \\ \hline  
1 & 25 & 10.2 & nc & - & nc & - & nc & - & 15 & 8.7 & nc & - & nc & - & nc & - \\ 
2 & 24 & 9.8 & nc & - & nc & - & nc & - & 15 & 8.6 & nc & - & nc & - & nc & - \\ 
3 & 24 & 9.9 & nc & - & nc & - & nc & - & 15 & 8.8 & nc & - & nc & - & nc & - \\ 
4 & 24 & 9.6 & nc & - & nc & - & nc & - & 15 & 8.3 & nc & - & nc & - & nc & - \\ 
6 & 23 & 9.8 & nc & - & nc & - & nc & - & 15 & 8.2 & nc & - & nc & - & nc & - \\ 
8 & 23 & 8.2 & nc & - & nc & - & nc & - & 15 & 7.7 & nc & - & nc & - & nc & - \\ 
12 & 23 & 9.8 & nc & - & nc & - & nc & - & 15 & 8.7 & nc & - & nc & - & nc & - \\ 
16 & 23 & 11.5 & 84 & 42.3 & nc & - & nc & - & 15 & 9.8 & 62 & 40.6 & nc & - & nc & - \\ 
24 & 22 & 14.2 & 53 & 34.3 & 91 & 59.0 & 94 & 61.6 & 14 & 11.2 & 38 & 30.4 & 54 & 43.3 & 57 & 45.8 \\ 
32 & 21 & 16.7 & 36 & 28.7 & 50 & 39.8 & 52 & 41.5 & 13 & 12.3 & 25 & 23.7 & 30 & 28.5 & 31 & 29.4 \\ 
40 & 20 & 19.0 & 33 & 31.3 & 40 & 38.0 & 41 & 39.0 & 13 & 14.3 & 22 & 24.3 & 24 & 26.6 & 25 & 27.7 \\ 
48 & 19 & 21.1 & 30 & 33.3 & 34 & 37.8 & 35 & 39.0 & 13 & 16.4 & 20 & 25.3 & 21 & 26.7 & 21 & 26.6 \\ 
56 & 18 & 23.0 & 26 & 33.2 & 29 & 36.9 & 30 & 38.3 & 12 & 17.2 & 18 & 25.7 & 19 & 27.0 & 19 & 27.1 \\ 
64 & 17 & 24.5 & 22 & 31.9 & 25 & 36.2 & 26 & 37.6 & 11 & 17.6 & 15 & 23.9 & 16 & 25.7 & 16 & 25.6 \\ 
\hline
&\multicolumn{16}{|c|}{Coarse grid, $20 \times 20$} \\ 
\hline  
1 & 24 & 6.8 & nc & - & nc & - & nc & - & 15 & 6.4 & nc & - & nc & - & nc & - \\ 
2 & 24 & 6.6 & nc & - & nc & - & nc & - & 15 & 6.2 & nc & - & nc & - & nc & - \\ 
3 & 23 & 6.7 & nc & - & nc & - & nc & - & 15 & 6.3 & nc & - & nc & - & nc & - \\ 
4 & 23 & 6.2 & nc & - & nc & - & nc & - & 15 & 6.2 & nc & - & nc & - & nc & - \\ 
6 & 23 & 7.1 & nc & - & nc & - & nc & - & 15 & 6.8 & 99 & 44.6 & nc & - & nc & - \\ 
8 & 22 & 7.5 & 86 & 28.9 & nc & - & nc & - & 15 & 7.4 & 62 & 30.0 & nc & - & nc & - \\ 
12 & 21 & 8.6 & 50 & 20.1 & 86 & 35.1 & 89 & 36.3 & 14 & 7.9 & 35 & 19.3 & 53 & 29.4 & 55 & 30.2 \\ 
16 & 19 & 9.2 & 34 & 16.2 & 49 & 23.7 & 51 & 24.6 & 13 & 8.3 & 23 & 14.4 & 29 & 18.4 & 31 & 19.3 \\ 
24 & 15 & 9.6 & 19 & 12.1 & 24 & 15.5 & 25 & 16.1 & 10 & 8.0 & 13 & 10.2 & 15 & 12.0 & 15 & 11.8 \\ 
32 & 14 & 11.6 & 16 & 13.2 & 18 & 14.8 & 18 & 15.0 & 9 & 9.0 & 11 & 10.7 & 11 & 10.9 & 11 & 10.7 \\ 
40 & 12 & 12.6 & 13 & 13.5 & 14 & 14.8 & 14 & 14.6 & 8 & 9.6 & 9 & 10.6 & 9 & 10.8 & 9 & 10.7 \\ 
48 & 11 & 14.1 & 12 & 15.3 & 12 & 15.3 & 12 & 15.2 & 7 & 9.9 & 8 & 11.2 & 8 & 11.3 & 7 & 9.9 \\ 
56 & 10 & 15.4 & 11 & 16.6 & 11 & 16.8 & 11 & 16.7 & 7 & 11.9 & 7 & 11.7 & 7 & 12.0 & 7 & 11.6 \\ 
64 & 9 & 16.3 & 10 & 17.9 & 10 & 18.2 & 9 & 16.3 & 6 & 11.9 & 7 & 13.8 & 6 & 11.9 & 6 & 11.8 \\ 
\hline
{\small mg$_1$} & 36 & 9.3 & nc & - & nc & - & nc & - & 16 & 4.2 & nc & - & nc & - & nc & - \\ 
{\small mg$_2$} & 34 & 14.5 & nc & - & nc & - & nc & - & 21 & 9.1 & nc & - & nc & - & nc & - \\ 
{\small mg$_3$} & 27 & 18.0 & nc & - & nc & - & nc & - & 18 & 12.0 & nc & - & nc & - & nc & - \\ 
\hline
\end{tabular}
\caption{Test 1. Average number of iterations $\bar{N}$ with time of solution tm$_s$}
\label{table-t1pr}
\end{table}

\begin{table}[h!]
\centering
\begin{tabular}{|c||cc|cc|cc|cc||cc|cc|cc|cc|}
\hline
\multirow{ 2}{*}{$J$} 
 & $\bar{N}$ & tm$_{s}$
 & $\bar{N}$ & tm$_{s}$
 & $\bar{N}$ & tm$_{s}$
 & $\bar{N}$ & tm$_{s}$
 & $\bar{N}$ & tm$_{s}$
 & $\bar{N}$ & tm$_{s}$
 & $\bar{N}$ & tm$_{s}$
 & $\bar{N}$ & tm$_{s}$\\ 
  & \multicolumn{2}{c|}{$10^3$} 
  & \multicolumn{2}{c|}{$10^6$}  
  & \multicolumn{2}{c|}{$10^9$}   
  & \multicolumn{2}{c||}{$10^{12}$}
  & \multicolumn{2}{c|}{$10^3$} 
  & \multicolumn{2}{c|}{$10^6$}  
  & \multicolumn{2}{c|}{$10^9$}   
  & \multicolumn{2}{c|}{$10^{12}$}\\ \hline
&\multicolumn{8}{|c||}{Jacobi smoother}&\multicolumn{8}{|c|}{Gauss-Seidel smoother} \\ 
\hline
& \multicolumn{16}{|c|}{Coarse grid, $10 \times 10$} \\ \hline  
1 & 26 & 7.5 & nc & - & nc & - & nc & - & 15 & 8.1 & nc & - & nc & - & nc & - \\ 
2 & 25 & 7.5 & nc & - & nc & - & nc & - & 15 & 8.1 & nc & - & nc & - & nc & - \\ 
3 & 24 & 7.4 & nc & - & nc & - & nc & - & 15 & 8.0 & nc & - & nc & - & nc & - \\ 
4 & 24 & 7.7 & nc & - & nc & - & nc & - & 15 & 7.9 & nc & - & nc & - & nc & - \\ 
6 & 23 & 7.8 & nc & - & nc & - & nc & - & 15 & 8.1 & nc & - & nc & - & nc & - \\ 
8 & 23 & 8.5 & nc & - & nc & - & nc & - & 15 & 7.7 & nc & - & nc & - & nc & - \\ 
12 & 23 & 9.9 & nc & - & nc & - & nc & - & 15 & 8.7 & 93 & 54.8 & nc & - & nc & - \\ 
16 & 23 & 11.4 & 83 & 41.4 & nc & - & nc & - & 15 & 9.8 & 57 & 37.5 & 65 & 42.9 & 64 & 42.4 \\ 
24 & 22 & 14.2 & 52 & 33.5 & 48 & 31.0 & 48 & 31.0 & 14 & 11.2 & 38 & 30.4 & 29 & 23.3 & 28 & 22.9 \\ 
32 & 21 & 16.7 & 37 & 29.4 & 30 & 23.8 & 29 & 23.1 & 13 & 12.3 & 25 & 23.7 & 18 & 17.1 & 18 & 17.1 \\ 
40 & 20 & 18.9 & 32 & 30.4 & 24 & 45.1 & 23 & 22.6 & 13 & 14.4 & 22 & 24.2 & 15 & 16.5 & 14 & 15.6 \\ 
48 & 19 & 21.0 & 28 & 31.0 & 20 & 22.1 & 20 & 22.3 & 12 & 15.1 & 19 & 24.0 & 12 & 15.1 & 12 & 15.2 \\ 
56 & 18 & 22.9 & 26 & 33.2 & 18 & 53.1 & 17 & 21.6 & 12 & 17.1 & 17 & 24.2 & 11 & 15.7 & 11 & 15.7 \\ 
64 & 17 & 24.5 & 23 & 33.3 & 16 & 23.1 & 16 & 23.1 & 11 & 17.6 & 16 & 25.7 & 10 & 16.0 & 10 & 16.0 \\ 
\hline
&\multicolumn{16}{|c|}{Coarse grid, $20 \times 20$} \\ 
\hline  
1 & 24 & 6.7 & nc & - & nc & - & nc & - & 16 & 6.7 & nc & - & nc & - & nc & - \\ 
2 & 24 & 6.8 & nc & - & nc & - & nc & - & 16 & 6.8 & nc & - & nc & - & nc & - \\ 
3 & 23 & 6.3 & nc & - & nc & - & nc & - & 15 & 6.8 & nc & - & nc & - & nc & - \\ 
4 & 23 & 6.2 & nc & - & nc & - & nc & - & 15 & 6.3 & nc & - & nc & - & nc & - \\ 
6 & 23 & 6.9 & nc & - & nc & - & nc & - & 15 & 6.8 & 86 & 39.0 & 98 & 44.9 & 98 & 44.1 \\ 
8 & 22 & 7.4 & 84 & 28.3 & nc & - & nc & - & 15 & 7.3 & 54 & 26.3 & 63 & 30.7 & 62 & 30.1 \\ 
12 & 21 & 8.6 & 49 & 19.9 & 45 & 18.5 & 45 & 18.0 & 14 & 7.8 & 34 & 19.0 & 27 & 15.1 & 27 & 14.8 \\ 
16 & 19 & 9.1 & 33 & 15.9 & 28 & 13.4 & 27 & 13.2 & 13 & 8.2 & 22 & 13.9 & 17 & 10.7 & 16 & 10.4 \\ 
24 & 16 & 10.2 & 21 & 13.5 & 14 & 8.9 & 14 & 8.9 & 11 & 8.7 & 14 & 11.1 & 9 & 7.1 & 9 & 7.1 \\ 
32 & 14 & 11.6 & 16 & 13.2 & 11 & 9.1 & 11 & 9.0 & 9 & 8.9 & 11 & 10.7 & 7 & 6.8 & 7 & 6.8 \\ 
40 & 12 & 12.6 & 13 & 13.4 & 9 & 9.4 & 9 & 9.3 & 8 & 9.5 & 9 & 10.7 & 6 & 7.2 & 6 & 7.1 \\ 
48 & 11 & 14.1 & 12 & 15.3 & 8 & 10.3 & 8 & 10.1 & 7 & 9.9 & 8 & 11.4 & 5 & 7.2 & 5 & 7.0 \\ 
56 & 10 & 15.3 & 11 & 16.9 & 7 & 10.8 & 7 & 10.6 & 7 & 11.8 & 7 & 11.7 & 5 & 8.5 & 4 & 6.7 \\ 
64 & 9 & 16.3 & 10 & 17.9 & 6 & 10.9 & 6 & 10.7 & 6 & 11.8 & 7 & 13.6 & 4 & 8.0 & 4 & 8.0 \\ 
 \hline
{\small mg$_1$} & 36 & 9.3 & nc & - & nc & - & nc & - & 17 & 4.4 & nc & - & nc & - & nc & - \\ 
{\small mg$_2$} & 31 & 13.3 & nc & - & nc & - & nc & - & 19 & 8.2 & nc & - & nc & - & nc & - \\ 
{\small mg$_3$} & 24 & 16.0 & nc & - & nc & - & nc & - & 16 & 10.6 & nc & - & nc & - & nc & - \\ 
\hline
\end{tabular}
\caption{Test 2. Average number of iterations $\bar{N}$ with time of solution tm$_s$}
\label{table-t2pr}
\end{table}

\begin{table}[h!]
\centering
\begin{tabular}{|c||cc|cc|cc|cc||cc|cc|cc|cc|}
\hline
\multirow{ 2}{*}{$J$} 
 & $\bar{N}$ & tm$_{s}$
 & $\bar{N}$ & tm$_{s}$
 & $\bar{N}$ & tm$_{s}$
 & $\bar{N}$ & tm$_{s}$
 & $\bar{N}$ & tm$_{s}$
 & $\bar{N}$ & tm$_{s}$
 & $\bar{N}$ & tm$_{s}$
 & $\bar{N}$ & tm$_{s}$\\ 
  & \multicolumn{2}{c|}{$10^3$} 
  & \multicolumn{2}{c|}{$10^6$}  
  & \multicolumn{2}{c|}{$10^9$}   
  & \multicolumn{2}{c||}{$10^{12}$}
  & \multicolumn{2}{c|}{$10^3$} 
  & \multicolumn{2}{c|}{$10^6$}  
  & \multicolumn{2}{c|}{$10^9$}   
  & \multicolumn{2}{c|}{$10^{12}$}\\ \hline
&\multicolumn{8}{|c||}{Jacobi smoother}&\multicolumn{8}{|c|}{Gauss-Seidel smoother} \\ 
\hline
& \multicolumn{16}{|c|}{Coarse grid, $10 \times 10$} \\ \hline  
1 & 19 & 7.2 & nc & - & nc & - & nc & - & 12 & 8.1 & 91 & 61.7 & 97 & 65.4 & 97 & 65.8 \\ 
2 & 19 & 7.5 & nc & - & nc & - & nc & - & 12 & 8.3 & 83 & 57.9 & 83 & 57.9 & 83 & 57.8 \\ 
3 & 19 & 7.8 & nc & - & nc & - & nc & - & 12 & 8.5 & 73 & 52.5 & 73 & 51.7 & 73 & 51.9 \\ 
4 & 19 & 8.2 & nc & - & nc & - & nc & - & 11 & 8.0 & 66 & 48.7 & 67 & 48.7 & 67 & 48.8 \\ 
6 & 18 & 8.2 & nc & - & nc & - & nc & - & 11 & 8.4 & 55 & 42.1 & 57 & 43.3 & 57 & 44.0 \\ 
8 & 26 & 9.3 & nc & - & nc & - & nc & - & 15 & 7.7 & 67 & 34.1 & 64 & 32.6 & 67 & 34.5 \\ 
12 & 25 & 10.7 & 96 & 41.5 & 99 & 42.3 & 99 & 42.5 & 15 & 8.7 & 48 & 27.9 & 49 & 28.5 & 49 & 28.7 \\ 
16 & 24 & 12.0 & 86 & 42.9 & 89 & 44.4 & 89 & 44.7 & 14 & 9.2 & 46 & 30.0 & 46 & 30.0 & 46 & 30.1 \\ 
24 & 23 & 14.9 & 62 & 40.0 & 62 & 40.1 & 62 & 40.2 & 13 & 10.4 & 32 & 25.5 & 32 & 25.5 & 32 & 25.6 \\ 
32 & 22 & 17.6 & 44 & 35.0 & 44 & 35.0 & 44 & 35.1 & 13 & 12.4 & 24 & 22.7 & 24 & 22.7 & 24 & 22.9 \\ 
40 & 21 & 20.0 & 40 & 37.9 & 40 & 37.9 & 40 & 38.0 & 12 & 13.3 & 21 & 23.1 & 21 & 23.1 & 21 & 23.2 \\ 
48 & 20 & 22.2 & 35 & 38.8 & 35 & 38.8 & 35 & 39.0 & 12 & 15.2 & 19 & 24.0 & 19 & 24.0 & 19 & 24.0 \\ 
56 & 19 & 24.3 & 31 & 39.5 & 32 & 40.7 & 32 & 40.8 & 11 & 15.7 & 17 & 24.3 & 17 & 24.3 & 17 & 24.3 \\ 
64 & 18 & 26.2 & 28 & 40.5 & 28 & 40.4 & 28 & 40.6 & 10 & 16.0 & 15 & 24.0 & 15 & 24.1 & 15 & 24.0 \\ 
\hline
&\multicolumn{16}{|c|}{Coarse grid, $20 \times 20$} \\ 
\hline  
1 & 26 & 5.6 & nc & - & nc & - & nc & - & 16 & 5.9 & nc & - & nc & - & nc & - \\ 
2 & 27 & 6.3 & nc & - & nc & - & nc & - & 16 & 6.2 & nc & - & nc & - & nc & - \\ 
3 & 26 & 6.4 & nc & - & nc & - & nc & - & 15 & 6.0 & 81 & 32.9 & 81 & 32.9 & 80 & 32.5 \\ 
4 & 26 & 6.9 & nc & - & nc & - & nc & - & 15 & 6.2 & 63 & 26.6 & 68 & 29.1 & 68 & 28.5 \\ 
6 & 24 & 7.1 & 94 & 28.4 & 94 & 28.7 & 94 & 29.2 & 14 & 6.3 & 43 & 19.7 & 43 & 19.7 & 43 & 19.7 \\ 
8 & 23 & 7.6 & 75 & 25.1 & 77 & 26.1 & 77 & 26.4 & 13 & 6.2 & 36 & 17.5 & 36 & 17.7 & 36 & 17.5 \\ 
12 & 20 & 7.9 & 47 & 18.9 & 47 & 19.0 & 47 & 19.2 & 12 & 6.5 & 25 & 14.0 & 25 & 14.0 & 25 & 14.0 \\ 
16 & 19 & 8.9 & 34 & 16.2 & 34 & 16.4 & 34 & 16.3 & 11 & 6.8 & 17 & 10.7 & 17 & 10.8 & 17 & 10.8 \\ 
24 & 16 & 10.1 & 21 & 13.5 & 21 & 13.4 & 21 & 13.5 & 10 & 7.8 & 12 & 9.6 & 12 & 9.5 & 12 & 9.5 \\ 
32 & 14 & 11.5 & 18 & 15.0 & 18 & 14.9 & 18 & 15.1 & 8 & 318.7 & 10 & 9.9 & 10 & 9.8 & 10 & 9.8 \\ 
40 & 12 & 12.4 & 14 & 14.7 & 14 & 14.5 & 14 & 14.7 & 7 & 8.1 & 8 & 9.5 & 8 & 9.5 & 8 & 9.5 \\ 
48 & 11 & 13.8 & 12 & 15.2 & 12 & 15.2 & 12 & 15.4 & 7 & 9.7 & 7 & 10.0 & 7 & 10.0 & 7 & 10.0 \\ 
56 & 10 & 15.2 & 11 & 16.9 & 11 & 16.8 & 11 & 16.8 & 6 & 9.9 & 6 & 9.9 & 6 & 10.2 & 6 & 10.1 \\ 
64 & 9 & 16.9 & 9 & 16.8 & 9 & 16.3 & 9 & 16.4 & 6 & 11.6 & 6 & 11.8 & 6 & 11.7 & 6 & 11.7 \\ 
\hline
{\small mg$_1$} & 31 & 8.1 & nc & - & nc & - & nc & - & 10 & 2.6 & 31 & 8.0 & 33 & 8.7 & 33 & 8.6 \\ 
{\small mg$_2$} & 30 & 12.8 & nc & - & nc & - & nc & - & 17 & 7.3 & nc & - & nc & - & nc & - \\ 
{\small mg$_3$} & 24 & 15.8 & nc & - & nc & - & nc & - & 14 & 9.3 & nc & - & nc & - & nc & - \\ 
\hline
\end{tabular}
\caption{Test 3. Average number of iterations $\bar{N}$ with time of solution tm$_s$}
\label{table-t3pr}
\end{table}

            
        

We consider the same test problems on a fine grid with $DOF_h=201,385$. We simulate with $t_{max} = 5 \cdot 10^{-6}$ using $N_t = 10$ time steps. For a coarse grid solver, we use the $10 \times 10$ and $20 \times 20$ multiscale approximation described in the previous section. We compare the performance of the proposed preconditioner with several classic multilevel AMG solvers from the PyAMG library. We consider three AMG preconditioners with pre- and post-smoothing to match the spectral two-level method, and default parameters from the PyAMG library (which have generally been tuned for robust performance): (1) mg$_1$ is a classical smoothed aggregation multilevel solver; (2) mg$_2$ is a classical AMG (Ruge-Stuben AMG) solver without second-pass coarsening; and (3) mg$_3$ is the same solver as mg$_2$ with second-pass coarsening for the C/F splitting (typically resulting in slower coarsening, but improved convergence). We investigate an average number of iterations for PCG per time iteration, $\bar{N} = \bar{N}_{tot}/N_t$, where $\bar{N}_{tot}$ is the total number of iterations. The results are presented with the time required for online computations (tm$_s$) for varying numbers of multiscale basis functions that directly affect the accuracy of the coarse grid approximation and size of the coarse grid system.

First and foremost, Tables \ref{table-t1pr}, \ref{table-t2pr}, and \ref{table-t3pr} demonstrate that the two-grid preconditioning built on a GMSFEM spectral coarse space is robust with respect to anisotropy. For all test cases, with sufficient coarse grid basis functions on a $20\times 20$ coarse grid, excellent convergence rates can be obtained that are independent of anisotropy. On a $10\times10$ coarse grid, there is modest growth in iteration count from $10^3$ to $10^6$ anistropy ratio, even for a large number of basis functions, but robust convergence is still obtained at anisotropy ratio $10^6$, which does not increase as the anisotropy ratio further increases. In general, the $10\times10$ coarse grid typically requires $2-3\times$ as many basis functions as the $20\times20$ coarse grid (usually closer to two) for comparable convergence. But, if basis functions can be computed offline, this would still be more efficient in practice, as the coarse grid operator is still smaller. We also point out that iteration counts monotonically decreases with increased coarse basis functions in almost all cases; in particular, indicating that convergence does not stall with respect to local basis functions, and can seemingly always be improved with more, down to $\mathcal{O}(1)$ iterations in many cases. The solve time also generally decreases with increased basis functions, at some point stagnating but rarely increasing notably. Again, if basis functions can be computed offline, this indicates that generally choosing more local basis functions will make for a rapidly converging and robust method.

Finally, we point out that Test 3 (see Table \ref{table-t3pr}), which appears to have the most complex field lines, generally observes the best and most robust convergence, particularly for large anisotropies. Test 3 also benefits significantly from symmetric Gauss-Seidel relaxation compared with Jacobi. Although iteration counts are fairly large, with Gauss Seidel we are able to converge for only 1 coarse basis function on the $10\times 10$ grid. We believe this is due to a large percentage of acyclic/open field lines in test 3, as opposed to the exclusively closed field lines in Tests 1 and 2. Closed field lines make the linear system particularly more ill-conditioned and challenging to solve, as in one time step heat will traverse the domain many times. In contrast, with open field lines there is a particular benefit to a Gauss-Seidel relaxation, which can relax on a given field line across the whole domain, boundary to boundary, analogous to an approximate line relaxation.

\section{Conclusion}

Here we present a multiscale finite element method to solve the heat flux problem with strong anisotropy. The presented method is based on constructing spectral multiscale basis functions and Galerkin coupling on the coarse grid. We have demonstrated that the presented basis functions are aligned with magnetic fields and produce excellent approximation for problems with high contrast between parallel and perpendicular flow directions. We demonstrate the potential of the multiscale basis functions in two contexts, (i) as a surrogate model with significantly reduced size, and (ii) as the coarse-grid correction in a two-level preconditioner. The latter is particularly relevant for extreme anisotropies as considered in this paper, because efficient implicit solvers remain a largely open question. The significantly reduced system size provided by the multiscale basis functions makes direct or sparse approximate inverses computationally tractable on the coarse grid, and we utilize only standard pointwise Jacobi or Gauss-Seidel relaxation on the fine grid. 

Numerical results demonstrate the method's convergence with respect to the number of spectral basis functions and coarse grid resolution. We have presented results for three test cases with different magnetic field distributions and varying anisotropy ratios.
The results show how a careful choice of the basis functions, with respect to the given magnetic field, can produce very good results for problems with very high anisotropy. Future work will consider more challenging and realistic test cases, applying the multiscale methodology to more specialized discretizations, e.g. \cite{wimmer2023fast,wimmer2024fast}, and developing efficient computational approaches for magnetic fields that evolve nonlinearly with the physical variables. 

\section*{Acknowledgements}
BSS and GW work was supported by the Laboratory Directed Research and Development program of Los Alamos National Laboratory under project number 20240261ER. LA-UR-24-32768.

\appendix
\section{Multiscale space}\label{app:ms}


We define $w^n \in V_H$ as the elliptic projection of $T_h^n \in V_h$  ($w^n = \Pi \ T_h^n$) that satisfies \eqref{pgest} and
\[
a(T_h^n - w^n, v) = 0, \quad \forall v \in V_H.
\]

For any $v \in V_H$, we have the following error equation
\[
m(T_h^n -T^n_{ms} , v) - m(T_h^{n-1} - T^{n-1}_{ms}, v) + \tau a(T_h^n - T^n_{ms}, v) 
= 
m(T_h^n  - T_h^{n-1}, v) + \tau a(T_h^n, v)  - \tau (f, v) = 0.
\]
We set $v = w^n - T^n_{ms} =  (T_h^n - T^n_{ms}) - (T_h^n - w^n)$ and obtain
\[
\begin{split}
&\underbrace{m(T_h^n -T^n_{ms}, T_h^n - T^n_{ms}) }_{I_1}
+ \tau 
\underbrace{a(T_h^n - T^n_{ms}, T_h^n - T^n_{ms}) }_{I_2}\\
&\quad \quad  = 
 \underbrace{m(T_h^n -T^n_{ms}, T_h^n - w^n) }_{I_3}
+ 
\underbrace{m(T_h^{n-1} - T^{n-1}_{ms},w^n - T^n_{ms}) }_{I_4}
+ \tau 
\underbrace{a(T_h^n - T^n_{ms}, T_h^n - w^n) }_{I_5}.
\end{split}
\]

For the first and second terms, we have
\[
I_1: \quad 
m(T_h^n -T^n_{ms}, T_h^n - T^n_{ms})
=||T_h^n - T^n_{ms}||_{M_h}^2, 
\]\[
I_2: \quad 
a(T_h^n - T^n_{ms}, T_h^n - T^n_{ms})
=||T_h^n - T^n_{ms}||_{A_h}^2.
\]

For the third and fifth terms by  Cauchy–Schwarz and Young's inequalities, we obtain
\[
I_3: \quad 
m(T_h^n -T^n_{ms},T_h^n - w^n)
 \leq
\frac{1}{4 \delta_1} ||T_h^n -T^n_{ms}||^2_{M_h} 
+ \delta_1 ||T_h^n - w^n||^2_{M_h},
\]\[
I_5: \quad 
a(T_h^n - T^n_{ms}, T_h^n - w^n)
\leq
\frac{1}{4 \delta_2}||T_h^n -T^n_{ms}||^2_{A_h} 
+ \delta_2 ||T_h^n - w^n||^2_{A_h}.
\]

For the forth term with $w^n - T^n_{ms} = (w^n - T_h^n) + (T_h^n - T^n_{ms})$, we have
\[
\begin{split}
I_4: \quad 
m(T_h^{n-1} - T^{n-1}_{ms}, & w^n - T^n_{ms})\\
&=
m(T_h^{n-1} - T^{n-1}_{ms}, w^n - T_h^n)
+
m(T_h^{n-1} - T^{n-1}_{ms}, T_h^n - T^n_{ms})\\
& \leq
\frac{1}{2 \delta_3}||T_h^{n-1} - T^{n-1}_{ms}||^2_{M_h} 
+ \delta_3 ||w^n - T_h^n||^2_{M_h} 
+ \delta_3 ||T_h^n - T^n_{ms}||^2_{M_h}.
\end{split}
\]

Then, we have 
\[
\begin{split}
&||T_h^n - T^n_{ms}||^2_{M_h} 
+ \tau ||T_h^n - T^n_{ms}||^2_{A_h}\\
&\quad \quad \leq 
\frac{1}{2 \delta_3} ||T_h^{n-1} - T^{n-1}_{ms}||^2_{M_h} 
+  \delta_3 ||w^n - T_h^n||^2_{M_h} 
+  \delta_3 ||T_h^n - T^n_{ms}||^2_{M_h}\\
&\quad \quad 
+ \frac{1}{4 \delta_1} ||T_h^n -T^n_{ms}||^2_{M_h} 
+ \delta_1 ||T_h^n - w^n||^2_{M_h}
+ \frac{\tau}{4 \delta_2}||T_h^n -T^n_{ms}||^2_{A_h} 
+ \tau \delta_2 ||T_h^n - w^n||^2_{A_h}\\
&\quad \quad =
\frac{1}{2 \delta_3} ||T_h^{n-1} - T^{n-1}_{ms}||^2_{M_h} 
+ \left( \delta_3 + \delta_1 \right) ||T_h^n - w^n||^2_{M_h} 
+ \left( \delta_3+ \frac{1}{4 \delta_1}\right) ||T_h^n - T^n_{ms}||^2_{M_h}\\
&\quad \quad 
+ \frac{\tau}{4 \delta_2}||T_h^n -T^n_{ms}||^2_{A_h} 
+ \tau \delta_2 ||T_h^n - w^n||^2_{A_h}.
\end{split}
\]

We can rewrite the estimate as follows
\[
\begin{split}
&  
\left( 1 -   \delta_3 - \frac{1}{4 \delta_1} \right) 
||T_h^n - T^n_{ms}||^2_{M_h} 
+ \tau 
\left( 1 - \frac{1}{4 \delta_2} \right) 
||T_h^n - T^n_{ms}||^2_{A_h}\\
&\quad \quad \leq 
\frac{1}{2 \delta_3} ||T_h^{n-1} - T^{n-1}_{ms}||^2_{M_h} 
+ \left(  \delta_3 + \delta_1 \right) ||T_h^n - w^n||^2_{M_h} 
+ \tau \delta_2 ||T_h^n - w^n||^2_{A_h}.
\end{split}
\]

For $\delta_1 = 1$, $\delta_2 = 1/2$ and $\delta_3 = 1/4$, we obtain
\[ 
||T_h^n - T^n_{ms}||^2_{M_h} 
+  \tau ||T_h^n - T^n_{ms}||^2_{A_h}
\preceq 
||T_h^{n-1} - T^{n-1}_{ms}||^2_{M_h}  
+ ||T_h^n - w^n||^2_{M_h} 
+  \tau||T_h^n - w^n||^2_{A_h}.
\]

Therefore
\[
||T_h^n - T^n_{ms}||^2_{M_h} 
+ \tau \sum_{k=1}^n ||T_h^n - T^n_{ms}||^2_{A_h}
\preceq
|| T_h^0 - T^0_{ms} ||^2_{M_h} 
+ \sum_{k=1}^n  \left(
||T_h^k - w^k||^2_{M_h}  + \tau ||T_h^k - w^k||^2_{A_h}
\right).
\]


Then, using estimates \eqref{pgest}
\[
||T_h^k - w^k||_{A_h}^2 \leq 
 \left( \frac{1}{H^2 \lambda_{J+1}^2} + \frac{1}{\lambda_{J+1}} \right) ||B_h T_h^k||^2_{D_h},
\]\[
||T_h^k - w^k||^2_{M_h}  \preceq 
||T_h^k - w^k||^2_{D_h} \leq 
\frac{1}{\lambda_{J+1}^2}  ||B_h T_h^k||^2_{D_h},
\]
we obtain the result
\[
\begin{split}
||T_h^n - T^n_{ms}||^2_{M_h} 
&+ \tau \sum_{k=1}^n ||T_h^n - T^n_{ms}||^2_{A_h}\\
&\preceq
|| T_h^0 - T^0_{ms} ||^2_{M_h} 
+  \sum_{k=1}^n \left( 
\frac{1}{\lambda_{J+1}^2} 
+ 
\tau \left(
 \frac{1}{H^2 \lambda_{J+1}^2} + \frac{1}{\lambda_{J+1}} \right)  \right) || B_h T_h^k||^2_{D_h} \\
& \leq
|| T_h^0 - T^0_{ms} ||^2_{M_h} 
+ \tau \sum_{k=1}^n   \frac{1}{\lambda_{J+1}}   || B_h T_h^k||^2_{D_h}. 
\end{split}
\]
By scaling eigenvalues with $H^{-2}$ ($\lambda_{J+1} = H^{-2} \Lambda^*$) for moving to size one domain \cite{abreu2019convergence}
\[
||T_h^n - T^n_{ms}||^2_{M_h} 
+ \tau \sum_{k=1}^n ||T_h^n - T^n_{ms}||^2_{A_h}
\preceq
|| T_h^0 - T^0_{ms} ||^2_{M_h}  \\
+ \tau  \sum_{k=1}^n \frac{H^2}{\Lambda^*}   || B_h T_h^k||^2_{D_h}
\]
Finally, under some additional regularity and appropriate initial conditions, we have an error of order $\mathcal{O}(\tau \frac{H^2}{\Lambda^*})$, where $\Lambda^*$ is responsible for covering a highly anisotropic flow by a spectral coarse space.

\section{Two-grid convergence}\label{app:two-grid}


Let $||v||_{Q_h}$ be a $\tau$-weighted $H_1$ norm \cite{akkutlu2015multiscale}
\[
||v||_{Q_h}^2 = v^T Q_h v 
= v^T \left( \frac{1}{\tau} M_h + A_h \right) v 
= \frac{1}{\tau} ||v||^2_{M_h} + ||v||^2_{A_h}.
\]

The two-grid error propagation can be expressed as follows \cite{falgout2005two, vassilevski2008multilevel}
\[
E_{TG} = I - C^{-1}_{TG} Q_h 
\]
with \[
C^{-1}_{TG} = \bar{S}^{-1} + (I - S^{-T} Q_h) P Q^{-1}_H P^T (I - Q_h S^{-1}),
\quad  
\bar{S} = S (S + S^T  - Q_h)^{-1} S^T.
\]

We have
\[
||v||^2_{Q_h} \leq v^T C_{TG} v \leq K_{TG} ||v||^2_{Q_h},
 \quad \textnormal{where }
 K_{TG} \coloneqq \text{Cond} (C_{TG}^{-1} Q_h),
\]
and
\[
0 \leq v^T Q_h E_{TG} v \leq \left( 1 - \frac{1}{K_{TG}} \right) ||v||^2_{Q_h}.
\] 

For typical smoothers (Jacobi, Gauss-Seidel), we have $Q_h = D_Q - N_Q - N_Q^T$, where $D_Q$ is the diagonal of $Q_h$ and $-N_Q$ is the strictly lower triangular part of $Q_h$. 
Then $S = D_Q - N_Q$ for Gauss-Seidel smoother and $S = D_Q$ for Jacobi and $\tilde{S}$ is the symmetric matrix that is spectrally equivalent to $D_Q$ \cite{vassilevski2008multilevel}. 
Moreover, we have spectral equivalency of smoother $S$ to $D_Q = \text{diag} (Q_h)$.

Based on the theory of two-grid method in \cite{vassilevski2008multilevel, brezina2011smoothed, falgout2005two, xu2002method}, we have
\[
K_{TG} = \text{sup}_v \frac{||v - \Pi v||^2_{D_Q}}{||v||^2_{Q_h}}.
\]
Using estimate \eqref{pgest}, we obtain
\[
||v - \Pi v ||_{D_Q}^2 
\leq  
\left(1 + \frac{C}{\tau} \right) ||v - \Pi v ||_{D_h}^2
\leq
\left(1 + \frac{C}{\tau} \right) \frac{H^2}{\Lambda^*} ||v||^2_{A_h}
\leq 
\left(1 + \frac{C}{\tau} \right) \frac{H^2}{\Lambda^*} ||v||^2_{Q_h}
\]
and we have $K_{TG} = \left(1 + \frac{C}{\tau} \right) \frac{H^2}{\Lambda^*}$.

\bibliographystyle{plain}
\bibliography{lit}

\end{document}